\documentclass[A4j,11pt]{article}
\usepackage{graphics}
\usepackage{amsmath}
\usepackage{amssymb}
\usepackage{latexsym}
\usepackage{amsfonts}
\begin{document}

\centerline{{\huge\bf A Cartan type identity for}}
\centerline{{\huge\bf isoparametric hypersurfaces}}
\centerline{{\huge\bf in symmetric spaces}}

\vspace{0.8truecm}

\centerline{\Large Naoyuki Koike}

\vspace{0.4truecm}

\begin{abstract}
In this paper, we obtain a Cartan type identity for curvature-adapted 
isoparametric hypersurfaces in symmetric spaces of compact type or 
non-compact type.  This identity is a generalization of 
Cartan-D'Atri's identity for curvature-adapted(=amenable) 
isoparametric hypersurfaces in rank one symmetric spaces.  
Furthermore, by using the Cartan type identity, we show that 
certain kind of curvature-adapted isoparametric hypersurfaces 
in a symmetric space of non-compact type are principal orbits of 
Hermann actions.  
\end{abstract}


\hspace{2.3truecm}complex focal radius, Hermann action


\vspace{0.4truecm}

\section{Introduction}
An isoparametric hypersurface in a (general) Riemannian manifold  is a 
connected hypersurface whose sufficiently close parallel 
hypersurfaces are of constant mean curvature (see [HLO] for example).  
In this paper, we assume that all isoparametric hypersurfaces are complete.  
It is known that all isoparametric hypersurfaces in a symmetric space of 
compact type are equifocal in the sense of [TT] and that, 
conversely all equifocal hypersurfaces are isoparametric (see [HLO]).  
Also, it is known that all isoparametric hypersurfaces in a symmetric space 
of non-compact type are complex equifocal in the sense of [Koi2] and that, 
conversely, all curvature-adapted complex equifocal hypersurfaces are 
isoparametric (see Theorem 15 of [Koi3]), where the curvature-adaptedness 
implies that, for a unit normal vector $v$, the (normal) Jacobi operator 
$R(\cdot,v)v$ preserves the tangent space invariantly and commutes with 
the shape operator $A$ for $v$, where $R$ is the curvature tensor of 
the ambient space.  It is known that principal orbits of a Hermann action 
(i.e., the action of a symmetric subgroup of $G$) of cohomogeneity one on a 
symmetric space $G/K$ of compact type are curvature-adapted and equifocal 
(see ([GT]).  Hence they are isoparametric hypersurfaces.  
On the other hand, we [Koi4,7] showed that 
the principal orbits of a Hermann action 
(i.e., the action of a (not necessarily compact) symmetric subgroup of $G$) 
of cohomogeneity one on a symmetric space $G/K$ of non-compact type are 
curvature-adapted and complex equifocal, and they have no 
focal point of non-Euclidean type on the ideal boundary of $G/K$.  
Hence they are isoparametric hypersurfaces.  

For an isoparametric hypersurface $M$ in a real space form $N$ of constant 
curvature $c$, it is known that the following Cartan's identity holds:
$$\sum_{\lambda\in{\rm Spec}A\setminus\{\lambda_0\}}
\frac{c+\lambda\lambda_0}{\lambda-\lambda_0}
\times m_{\lambda}=0\leqno{(1.1)}$$
for any $\lambda_0\in{\rm Spec}A$, where $A$ is the shape operator of $M$ 
and ${\rm Spec}A$ is the spectrum of $A$, $m_{\lambda}$ is the multiplicity 
of $\lambda$.  
Here we note that all hypersurfaces in a real space form are 
curvature-adapted.  
In general cases, this identity is shown in algebraic method.  Also, 
It is shown in geometrical method in the following three cases:

\vspace{0.2truecm}

(i) $c=0,\,\,\lambda_0\not=0$,

(ii) $c>0,\,\,\lambda_0\,:\,$ any eigenvalue of $A_v$,

(iii) $c<0,\,\,\vert\lambda_0\vert>\sqrt{-c}$.  

\vspace{0.2truecm}

\noindent
In detail, it is shown by showing the minimality of the focal submanifold for 
$\lambda_0$ and using this fact.  

Let $H\curvearrowright G/K$ be a cohomogeneity one action of a compact group 
$H\,(\subset G)$ on a rank one symmetric space $G/K$ and $M$ a principal 
orbit of this action.  Since the $H$-action is of cohomogeneity one, it 
is hyperpolar.  Hence $M$ is an equifocal (hence isoparametric) hypersurface 
(see [HPTT]).  
In 1979, J. E. D'Atri [D] obtained a Cartan type identity for $M$ 
in the case where $M$ is amenable (i.e., curvature-adapted).  
On the other hand, in 1989-1991, J. Berndt [B1,2] obtained a Cartan type 
identity (in algebraic method) for curvature-adapted hypersurfaces with 
constant principal curvature in rank one symmetric spaces other than spheres 
and hyperbolic spaces.  
Here we note that, for a curvature-adapted hypersurface in a rank one 
symmetric space of non-compact type, it has constant principal curvature 
if and only if it is isoparametric.  

In this paper, we obtain the Cartan type identities 
for curvature-adapted isoparametric hypersurfaces in symmetric spaces and, 
furthermore, by using the Cartan type identity, we prove that 
certain kind of curvature-adapted isoparametric hypersurfaces 
in a symmetric space of non-compact type are principal orbits of 
Hermann actions.  
Let $M$ be a hypersurface in a symmetric space $N=G/K$ of 
compact type or non-compact type and $v$ a unit normal vector field of $M$.  
Set $R(v_x):=R(\cdot,v_x)v_x\vert_{T_xM}$, where $R$ is the curvature tensor 
of $N$.  For each $r\in \Bbb R$, we define a function $\tau_r$ over 
$[0,\infty)$ by 
$$\tau_r(s):=
\left\{
\begin{array}{cc}
\displaystyle{\frac{\sqrt s}{\tan(r\sqrt s)}} & \displaystyle{(s>0)}\\
\displaystyle{\frac 1r} & \displaystyle{(s=0)}
\end{array}\right.$$
Also, for each $r\in\Bbb C$, we define a complex-valued function 
$\hat{\tau}_r$ over $(-\infty,0]$ by 
$$\hat{\tau}_r(s):=
\left\{
\begin{array}{cc}
\displaystyle{\frac{{\bf i}\sqrt{-s}}{\tan({\bf i}r\sqrt{-s})}} & 
\displaystyle{(s<0)}\\
\displaystyle{\frac 1r} & \displaystyle{(s=0),}
\end{array}\right.$$
where ${\bf i}$ is the imaginary unit.  
First we prove the following Cartan type identity 
for a curvature-adapted isoparametric hypersurface in 
a simply connected symmetric space of compact type.  

\vspace{0.3truecm}

\noindent
{\bf Theorem A.} {\sl Let $M$ be a curvature-adapted isoparametric 
hypersurface in a simply connected symmetric space $N:=G/K$ of compact type.  
For each focal radius $r_0$ of $M$, we have 
$$\sum_{(\lambda,\mu)\in S_{r_0}^x}
\frac{\mu+\lambda \tau_{r_0}(\mu)}{\lambda-\tau_{r_0}(\mu)}\times 
m_{\lambda,\mu}=0,\leqno{(1.2)}$$
where $S_{r_0}^x:=\{(\lambda,\mu)\in {\rm Spec}A_x\times{\rm Spec}R(v_x)\,
\vert\,{\rm Ker}(A_x-\lambda I)\cap{\rm Ker}(R(v_x)-\mu I)\not=\{0\},\,\,
\lambda\not=\tau_{r_0}(\mu)\}$ and $m_{\lambda,\mu}:={\rm dim}
({\rm Ker}(A_x-\lambda I)\cap{\rm Ker}(R(v_x)-\mu I))$.}

\vspace{0.3truecm}

\noindent
{\bf Remark 1.1.} 
(i) If ${\rm Ker}(A_x-\lambda_0 I)\cap{\rm Ker}(R(v_x)-\mu_0 I)$ is included 
by the focal space for the focal radius $r_0$, then we have 
$\tau_{r_0}(\mu_0)=\lambda_0$.  

(ii) If $G/K$ is a sphere of constant curvature $c$, then 
${\rm Spec}R(v_x)=\{c\}$ and $\tau_{r_0}(c)$ is equal to the principal 
curvature corresponding to $r_0$.  Hence the identity 
$(1.2)$ coincides with $(1.1)$.  

(iii) In the case where $G/K$ is a rank one symmetric 
space of compact type, the identity $(1.2)$ 
coincides with the identity obtained by J. E. D'Atri [D] (see Theorems 3.7 
and 3.9 of [D]).  

(iv) In the case where $G/K$ is a rank one symmetric 
space of compact type other than spheres, the identity $(1.2)$ 
is different from the identity obtained by J. Berndt [B1,2].  

\vspace{0.3truecm}

Next, in this paper, we prove the following Cartan type identity for 
a curvature-adapted isoparametric $C^{\omega}$-hypersurface in 
a symmetric space of non-compact type, where $C^{\omega}$ means the real 
analyticity.  

\vspace{0.3truecm}

\noindent
{\bf Theorem B.} {\sl Let $M$ be a curvature-adapted isoparametric 
$C^{\omega}$-hypersurface in a symmetric space $N:=G/K$ of non-compact type.  
Assume that 
$M$ has no focal point of non-Euclidean type on the ideal 
boundary $N(\infty)$ of $N$.  
Then $M$ admits a complex focal radius and , for each complex focal radius 
$r_0$ of $M$, we have 
$$\sum_{(\lambda,\mu)\in S_{r_0}^x}
\frac{\mu+\lambda\hat{\tau}_{r_0}(\mu)}
{\lambda-\hat{\tau}_{r_0}(\mu)}\times m_{\lambda,\mu}=0,\leqno{(1.3)}$$
where $S_{r_0}^x:=\{(\lambda,\mu)\in {\rm Spec}A_x\times{\rm Spec}R(v_x)\,
\vert\,{\rm Ker}(A_x-\lambda I)\cap{\rm Ker}(R(v_x)-\mu I)\not=\{0\},\,\,
\lambda\not=\hat{\tau}_{r_0}(\mu)\}$ and 
$m_{\lambda,\mu}:={\rm dim}({\rm Ker}(A_x-\lambda I)\cap
{\rm Ker}(R(v_x)-\mu I))$.}

\vspace{0.3truecm}

\noindent
{\bf Remark 1.2.} 
(i) The notion of a complex focal radius was introduced in [Koi2].  
This quantity indicates the position of a focal point of the complexification 
$M^{\bf c}\,(\subset G^{\bf c}/K^{\bf c})$ of a submanifold $M$ in a symmetric 
space $G/K$ of non-compact type (see [Koi3]).  

(ii) If ${\rm Ker}(A_x-\lambda_0 I)\cap{\rm Ker}(R(v_x)-\mu_0 I)$ is included 
by the focal space for the complex focal radius $r_0$, then we have 
$\hat{\tau}_{r_0}(\mu_0)=\lambda_0$.  

(iii) If $G/K$ is a hyperbolic space of constant curvature $c$, then 
${\rm Spec}R(v_x)=\{c\}$ and $\hat{\tau}_{r_0}(c)$ is equal to the principal 
curvature corresponding to $r_0$.  Hence the identity $(1.3)$ coincides with 
$(1.1)$.  

(iv) In the case where $G/K$ is a rank one symmetric 
space of non-compact type and $r_0$ is a real focal radius, the identity 
$(1.3)$ coincides with the identity obtained by J. E. D'Atri [D] 
(see Theorems 3.7 and 3.9 of [D]).  

(v) In the case where $G/K$ is a rank one symmetric 
space of non-compact type other than hyperbolic spaces, the identity $(1.3)$ 
is different from the identity obtained by J. Berndt [B1,2].  

(vi) For a curvature-adapted and isoparametric hypersurface $M$ in $G/K$, the 
following conditions ${\rm (a)}\sim{\rm (c)}$ are equivalent:

(a) $M$ has no focal point of non-Euclidean type on $N(\infty)$,  

(b) $M$ is proper complex equifocal in the sense of [Koi4],

(c) ${\rm Ker}(A_x\pm\sqrt{-\mu}I)\cap{\rm Ker}(R(v_x)-\mu I)=\{0\}$ holds 
for each $\mu\in{\rm Spec}R(v_x)\setminus\{0\}$.  

(vii) Principal orbits of a Hermann type action of cohomogeneity one on 
$G/K$ are curvature-adapted isoparametric $C^{\omega}$-hypersurface 
having no focal point of non-Euclidean type on $N(\infty)$ 
(see Theorem B of [Koi4] and the above (iii)).  

\vspace{0.3truecm}

The proof of Theorem B is performed by showing {\bf the minimality of the 
focal submanifold} $F:=\{\exp^{\perp}(({\rm Re}\,r_0)v_x+({\rm Im}\,r_0)Jv_x)
\,\vert\,x\in M^{\bf c}\}$ of the complexification $M^{\bf c}$ of 
$M$ (see Fig.1), where $\exp^{\perp}$ is the normal exponential map of 
the submanifold $M^{\bf c}$ in $G^{\bf c}/K^{\bf c}$, $J$ is the complex 
structure of $G^{\bf c}/K^{\bf c}$ and $v$ is a unit normal vector field 
of $M$ (in $G/K$).  
Here we note that $\exp^{\perp}(({\rm Re}\,r_0)v_x+({\rm Im}\,r_0)Jv_x)$ 
is equal to the point $\gamma_{v_x}^{\bf c}(r_0)$ of the complexified geodesic 
$\gamma_{v_x}^{\bf c}$ in $G^{\bf c}/K^{\bf c}$.  
In the case where $G/K$ is of rank greater than one and $M$ is not 
homogeneous, the proof of the minimality of $F$ is performed by showing 
{\bf the minimality of the lift} $\widetilde F:=(\pi\circ\phi)^{-1}(F)$ 
{\bf of} $F$ {\bf to the path space} $H^0([0,1],\mathfrak g^{\bf c})$, where 
$\phi$ is the parallel transport map for $G^{\bf c}$ (which is an 
anti-Kaehlerian submersion o $H^0([0,1],\mathfrak g^{\bf c})$ onto 
$G^{\bf c}$) and $\pi$ is the natural projection of $G^{\bf c}$ onto 
$G^{\bf c}/K^{\bf c}$ (which also is an anti-Kaehlerian submersion).  
Here we note that 
the minimality of $F$ is trivial in the case where $M$ is homogeneous.  
By using Theorem B, we prove the following fact for the number of 
distinct principal curvatures of a curvature-adapted isoparametric 
$C^{\omega}$-hypersurfaces in a symmetric sapce of 
non-compact type.  

\vspace{0.3truecm}

\centerline{
\unitlength 0.1in
\begin{picture}( 39.8500, 18.9300)( -0.5300,-23.0300)
%
\special{pn 8}%
\special{pa 1566 1968}%
\special{pa 2610 2146}%
\special{fp}%
%
\special{pn 8}%
\special{ar 2136 1956 986 220  5.7136066 6.2831853}%
\special{ar 2136 1956 986 220  0.0000000 3.6551879}%
%
\special{pn 8}%
\special{ar 2136 1956 978 214  4.4093875 4.4295387}%
\special{ar 2136 1956 978 214  4.4899921 4.5101432}%
\special{ar 2136 1956 978 214  4.5705966 4.5907477}%
\special{ar 2136 1956 978 214  4.6512011 4.6713523}%
\special{ar 2136 1956 978 214  4.7318057 4.7519568}%
\special{ar 2136 1956 978 214  4.8124102 4.8325613}%
\special{ar 2136 1956 978 214  4.8930147 4.9131659}%
\special{ar 2136 1956 978 214  4.9736193 4.9937704}%
\special{ar 2136 1956 978 214  5.0542238 5.0743749}%
\special{ar 2136 1956 978 214  5.1348283 5.1549795}%
\special{ar 2136 1956 978 214  5.2154329 5.2355840}%
\special{ar 2136 1956 978 214  5.2960374 5.3161885}%
\special{ar 2136 1956 978 214  5.3766419 5.3967931}%
\special{ar 2136 1956 978 214  5.4572465 5.4773976}%
\special{ar 2136 1956 978 214  5.5378510 5.5580021}%
%
\special{pn 8}%
\special{ar 2160 1948 978 214  3.7917162 3.8118673}%
\special{ar 2160 1948 978 214  3.8723207 3.8924718}%
\special{ar 2160 1948 978 214  3.9529253 3.9730764}%
\special{ar 2160 1948 978 214  4.0335298 4.0536809}%
%
\special{pn 20}%
\special{ar 1972 2290 452 790  3.3497801 3.8592609}%
%
\special{pn 8}%
\special{ar 1916 918 858 1344  2.0186185 3.0192679}%
%
\special{pn 8}%
\special{pa 1626 1752}%
\special{pa 1606 1726}%
\special{pa 1586 1700}%
\special{pa 1568 1674}%
\special{pa 1552 1648}%
\special{pa 1534 1620}%
\special{pa 1518 1594}%
\special{pa 1502 1564}%
\special{pa 1488 1536}%
\special{pa 1474 1508}%
\special{pa 1462 1478}%
\special{pa 1448 1448}%
\special{pa 1436 1420}%
\special{pa 1426 1390}%
\special{pa 1418 1358}%
\special{pa 1408 1328}%
\special{pa 1400 1296}%
\special{pa 1394 1266}%
\special{pa 1388 1234}%
\special{pa 1382 1202}%
\special{pa 1378 1170}%
\special{pa 1374 1140}%
\special{pa 1372 1108}%
\special{pa 1370 1076}%
\special{pa 1372 1044}%
\special{pa 1372 1012}%
\special{pa 1374 980}%
\special{pa 1376 948}%
\special{pa 1380 916}%
\special{pa 1386 884}%
\special{pa 1390 852}%
\special{pa 1396 828}%
\special{sp}%
%
\special{pn 8}%
\special{pa 1090 1070}%
\special{pa 1108 1044}%
\special{pa 1128 1018}%
\special{pa 1150 996}%
\special{pa 1170 972}%
\special{pa 1194 952}%
\special{pa 1218 930}%
\special{pa 1244 912}%
\special{pa 1270 892}%
\special{pa 1298 876}%
\special{pa 1326 860}%
\special{pa 1354 846}%
\special{pa 1382 830}%
\special{pa 1388 828}%
\special{sp}%
%
\special{pn 8}%
\special{ar 2114 1948 1008 214  4.1931684 4.3555232}%
%
\special{pn 20}%
\special{sh 1}%
\special{ar 1566 1968 10 10 0  6.28318530717959E+0000}%
\special{sh 1}%
\special{ar 1566 1968 10 10 0  6.28318530717959E+0000}%
%
\special{pn 13}%
\special{pa 1574 1968}%
\special{pa 2076 2054}%
\special{fp}%
\special{sh 1}%
\special{pa 2076 2054}%
\special{pa 2014 2024}%
\special{pa 2024 2046}%
\special{pa 2008 2062}%
\special{pa 2076 2054}%
\special{fp}%
%
\special{pn 13}%
\special{pa 1566 1962}%
\special{pa 1864 1622}%
\special{fp}%
\special{sh 1}%
\special{pa 1864 1622}%
\special{pa 1806 1658}%
\special{pa 1830 1662}%
\special{pa 1836 1684}%
\special{pa 1864 1622}%
\special{fp}%
%
\special{pn 8}%
\special{pa 2196 2290}%
\special{pa 2318 2098}%
\special{da 0.070}%
\special{sh 1}%
\special{pa 2318 2098}%
\special{pa 2266 2144}%
\special{pa 2290 2142}%
\special{pa 2298 2164}%
\special{pa 2318 2098}%
\special{fp}%
%
\special{pn 8}%
\special{pa 1754 2268}%
\special{pa 1550 2076}%
\special{da 0.070}%
\special{sh 1}%
\special{pa 1550 2076}%
\special{pa 1584 2136}%
\special{pa 1588 2112}%
\special{pa 1612 2106}%
\special{pa 1550 2076}%
\special{fp}%
\put(30.8300,-17.9100){\makebox(0,0)[lb]{$G/K$}}%
\put(21.2100,-23.0300){\makebox(0,0)[lt]{$\gamma_{v_x}$}}%
\put(18.2800,-22.8900){\makebox(0,0)[rt]{$M$}}%
\put(20.0900,-20.1100){\makebox(0,0)[lb]{$v_x$}}%
\put(18.0600,-16.3400){\makebox(0,0)[rb]{$Jv_x$}}%
\put(13.4500,-10.5600){\makebox(0,0)[rt]{$M^{\bf c}$}}%
\put(34.0000,-5.4600){\makebox(0,0)[lt]{$\gamma_{v_x}^{\bf c}$}}%
\put(15.4300,-19.1900){\makebox(0,0)[rt]{$x$}}%
\put(39.3200,-11.7200){\makebox(0,0)[lt]{in $G^{\bf c}/K^{\bf c}$}}%
\put(23.7700,-7.1600){\makebox(0,0)[rb]{$\gamma_{v_x}^{\bf c}(r_0)$}}%
%
\special{pn 13}%
\special{pa 1566 1968}%
\special{pa 2044 1692}%
\special{fp}%
\special{sh 1}%
\special{pa 2044 1692}%
\special{pa 1976 1708}%
\special{pa 1998 1718}%
\special{pa 1996 1742}%
\special{pa 2044 1692}%
\special{fp}%
%
\special{pn 8}%
\special{pa 1574 1968}%
\special{pa 2602 2140}%
\special{pa 3888 496}%
\special{pa 2806 410}%
\special{pa 2806 410}%
\special{pa 1574 1968}%
\special{fp}%
%
\special{pn 8}%
\special{pa 1580 1962}%
\special{pa 2918 1172}%
\special{da 0.070}%
%
\special{pn 20}%
\special{sh 1}%
\special{ar 2934 1164 10 10 0  6.28318530717959E+0000}%
\special{sh 1}%
\special{ar 2934 1164 10 10 0  6.28318530717959E+0000}%
%
\special{pn 8}%
\special{ar 3692 1008 790 584  2.8665774 2.9540410}%
\special{ar 3692 1008 790 584  3.0065191 3.0939827}%
\special{ar 3692 1008 790 584  3.1464608 3.2339244}%
\special{ar 3692 1008 790 584  3.2864025 3.3738661}%
%
\special{pn 20}%
\special{ar 3692 1008 782 612  1.9790358 2.8738906}%
%
\special{pn 8}%
\special{pa 3098 1798}%
\special{pa 2978 1956}%
\special{fp}%
\special{sh 1}%
\special{pa 2978 1956}%
\special{pa 3034 1914}%
\special{pa 3010 1914}%
\special{pa 3004 1890}%
\special{pa 2978 1956}%
\special{fp}%
\put(34.2900,-15.2800){\makebox(0,0)[lt]{$F$}}%
%
\special{pn 8}%
\special{pa 1692 1108}%
\special{pa 1980 1728}%
\special{da 0.070}%
\special{sh 1}%
\special{pa 1980 1728}%
\special{pa 1970 1658}%
\special{pa 1958 1680}%
\special{pa 1934 1676}%
\special{pa 1980 1728}%
\special{fp}%
\put(15.8000,-10.7900){\makebox(0,0)[lb]{$\frac{r_0}{\vert r_0\vert}v_x$}}%
%
\special{pn 8}%
\special{pa 2280 738}%
\special{pa 2926 1164}%
\special{da 0.070}%
\special{sh 1}%
\special{pa 2926 1164}%
\special{pa 2880 1112}%
\special{pa 2882 1136}%
\special{pa 2858 1144}%
\special{pa 2926 1164}%
\special{fp}%
%
\special{pn 8}%
\special{ar 1776 1208 414 876  2.1149712 2.2080666}%
\special{ar 1776 1208 414 876  2.2639239 2.3570193}%
\special{ar 1776 1208 414 876  2.4128765 2.5059720}%
\special{ar 1776 1208 414 876  2.5618292 2.5736661}%
%
\special{pn 8}%
\special{pa 1430 1692}%
\special{pa 2934 1164}%
\special{da 0.070}%
%
\special{pn 8}%
\special{pa 1484 1826}%
\special{pa 2934 1164}%
\special{da 0.070}%
\end{picture}%
\hspace{3truecm}}

\vspace{0.6truecm}

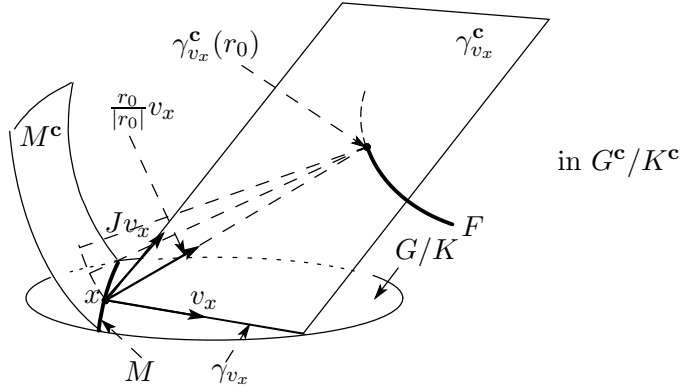
\centerline{{\bf Fig. 1.}}


\vspace{0.5truecm}

By using Theorem B, we prove the following main result.  

\vspace{0.5truecm}

\noindent
{\bf Theorem C.} {\sl Let $M$ be a curvature-adapted isoparametric 
$C^{\omega}$-hypersurface in a symmetric space $N$ of non-compact type.  
Assume that $M$ has no focal point of non-Euclidean type on $N(\infty)$.  
Then $M$ is a principal orbit of a Hermann action.}

\vspace{0.5truecm}

\noindent
{\bf Remark 1.3.} 
In this theorem, are indispensable both the condition of the 
curvature-adaptedness and the condition for the non-existenceness of 
non-Euclidean type focal point on the ideal boundary.  
In fact, we have the following examples.  
Let $G/K$ be an irreducible symmetric space of non-compact type such that 
the (restricted) root system of $G/K$ is non-reduced.  
Let $\mathfrak g=\mathfrak k+\mathfrak p$ ($\mathfrak g={\rm Lie}\,G,\,
\mathfrak k={\rm Lie}\,K$) be the Cartan decomposition associated with 
a symmetric pair $(G,K)$ and $\mathfrak a$ a maximal abelian subspace of 
$\mathfrak p$.  Also, let $\triangle_+$ be the positive root system of $G/K$ 
with respect to $\mathfrak a$ and $\Pi$ the simple root system of 
$\triangle_+$, where we fix a lexicographic ordering of the dual space 
$\mathfrak a^{\ast}$ of $\mathfrak a$.  
Set $\mathfrak n:=\sum_{\lambda\in\triangle_+}\mathfrak g_{\lambda}$ and 
$N:=\exp\,\mathfrak n$, where $\mathfrak g_{\lambda}$ is the root space for 
$\lambda$ and $\exp$ is the exponential map of $G$.  
If $G/K$ is of rank one, then any orbit of the $N$-action on $G/K$ is a full 
irreducible curvature-adapted isoparametric $C^{\omega}$-hypersurface but 
it has a focal point of non-Euclidean type on $N(\infty)$ (see [Koi9]).  
On the other hand, it is a principal orbit of no Hermann action.  
Thus, in this theorem, is indispensable the condition for the non-existenceness of a focal point of non-Euclidean type on the ideal boundary.  
Let $H_{\lambda}$ be the element of $\mathfrak a$ defined by 
$\langle H_{\lambda},\bullet\rangle=\lambda(\bullet)$.  
Assume that the (restricted) root system of $G/K$ is of type $(BC_n)$.  
Take an element $\lambda$ of $\Pi$ such that 
$2\lambda$ belongs to $\triangle_+$, and one-dimensional subspaces ${\it l}$ 
of ${\Bbb R}H_{\lambda}+\mathfrak g_{\lambda}$.  
Set $S:=\exp((\mathfrak a+\mathfrak n)\ominus{\it l})$, where 
$\exp$ is the exponential map of $G$ and $(\mathfrak a+\mathfrak n)\ominus{\it l}$ is 
the orthogonal complement of ${\it l}$ in $\mathfrak a+\mathfrak n$.  
Then $S$ is a subgroup of $AN:=\exp(\mathfrak a+\mathfrak n)$ and any orbit of 
the $S$-action on $G/K$ is a full irreducible isoparametric 
$C^{\omega}$-hypersurface but it is not curvature-adapted (see [Koi9]).  
Furthermore, we can find an orbit having no focal point of non-Euclidean type 
on $N(\infty)$ among orbits of the $S$-action.  
On the other hand, it is a principal orbit of no Hermann action.  
Thus the condition of the curvature-adaptedness is indispensable in this 
theorem.  

\vspace{0.5truecm}

In Section 2, we recall basic notions.  
In Section 3, we prove Theorem A.  
In Section 4, we define the mean curvature of a proper 
anti-Kaehlerian Fredholm submanifold and prepare a lemma to prove Theorem B.  
In Section 5, we prove Theorems B and C.  

\section{Basic notions}
In this section, we recall basic notions which are used in the proof 
of Theorems A and B.  
First we recall the notion of an equifocal hypersurface in a symmetric space.  
Let $M$ be a complete (oriented embedded) hypersurface in a symmetric space 
$N=G/K$ and fix a global unit normal vector field $v$ of $M$.  
Let $\gamma_{v_x}$ be the normal geodesic of $M$ with $\gamma_{v_x}'(0)=v_x$, 
where $x\in M$ and $\gamma_{v_x}'(0)$ is the velocity vector of $\gamma_{v_x}$ 
at $0$.  If $\gamma_{v_x}(s_0)$ is a focal point of $M$ along $\gamma_{v_x}$, 
then $s_0$ is called a {\it focal radius of} $M$ {\it at} $x$.  
Denote by ${\cal FR}_{M,x}$ the set of all focal radii of $M$ at $x$.  
If $M$ is compact and if ${\cal FR}_{M,x}$ is independent of the choice of 
$x$, then it is called an {\it equifocal hypersurface}.  
This notion is the hypersurface version of an equifocal submanifold defined 
in [TT].  

Next we recall the notion of a complex equifocal hypersurface in a symmetric 
space of non-compact type.  
Let $M$ be a complete (oriented embedded) hypersurface in a symmetric space 
$N=G/K$ of non-compact type and fix a global unit normal vector field $v$ of 
$M$.  Let $\mathfrak g$ be the Lie algebra of $G$ and $\theta$ be the Cartan 
involution of $G$ with ${\rm Fix}\,\theta=K$, where 
${\rm Fix}\,\theta$ is the fixed point group of $\theta$.  
Denote by the same symbol $\theta$ the involution of $\mathfrak g$ induced 
from $\theta$.  Set $\mathfrak p:={\rm Ker}(\theta+{\rm id})$.  
The subspace $\mathfrak p$ is identified with the tangent space 
$T_{eK}N$ of $N$ at $eK$, where $e$ is the identity element of $G$.  
Let $M$ be a complete (oriented embedded) hypersurface in $N$.  
Fix a global unit normal vector field $v$ of $M$.  
Denote by $A$ the shape operator of $M$ (for $v$).  
Take $X\in T_xM$ ($x=gK$).  
The 
$M$-Jacobi field $Y$ along 
$\gamma_x$ with $Y(0)=X$ (hence $Y'(0)=-A_xX$) is given by 
$$Y(s)=(P_{\gamma_x\vert_{[0,s]}}\circ(D^{co}_{sv_x}-sD^{si}_{sv_x}\circ 
A_x))(X),$$
where $P_{\gamma_x\vert_{[0,s]}}$ is the parallel translation along 
$\gamma_x\vert_{[0,s]}$, $D^{co}_{sv_x}$ (resp. $D^{si}_{sv_x}$) is given by 
$$\begin{array}{l}
\displaystyle{D^{co}_{sv_x}=
g_{\ast}\circ\cos({\bf i}{\rm ad}(sg_{\ast}^{-1}v_x))\circ g_{\ast}^{-1}}\\
\displaystyle{\left({\rm resp.}\,\,
D^{si}_{sv_x}=g_{\ast}\circ
\frac{\sin({\bf i}{\rm ad}(sg_{\ast}^{-1}v_x))}
{{\bf i}{\rm ad}(sg_{\ast}^{-1}v_x)}\circ g_{\ast}^{-1}\right).}
\end{array}$$ 
Here ${\rm ad}$ is the adjoint 
representation of the Lie algebra $\mathfrak g$ of $G$.  
All focal radii of $M$ at $x$ are catched 
as real numbers $s_0$ 
with ${\rm Ker}(D^{co}_{s_0v_x}-s_0D^{si}_{s_0v_x}\circ A_x)\not=\{0\}$.  
So, we [Koi2] defined the notion of a {\it complex focal radius of} $M$ 
{\it at} $x$ as a complex number $z_0$ with ${\rm Ker}(D^{co}_{z_0v_x}-
z_0D^{si}_{z_0v_x}\circ A_x^{{\bf c}})\not=\{0\}$, 
where $D^{co}_{z_0v_x}$ (resp. $D^{si}_{z_0v_x}$) 
is a ${\bf C}$-linear transformation of $(T_xN)^{\bf c}$ defined by 
$$\begin{array}{c}
\displaystyle{
D^{co}_{z_0v_x}=g^{\bf c}_{\ast}\circ\cos({\bf i}{\rm ad}^{\bf c}
(z_0g_{\ast}^{-1}v_x))\circ (g^{\bf c}_{\ast})^{-1}}\\
\displaystyle{\left({\rm resp.}\,\,\,\,
D^{si}_{sv_x}=g^{\bf c}_{\ast}\circ
\frac{\sin({\bf i}{\rm ad}^{\bf c}(z_0g_{\ast}^{-1}v_x))}
{{\bf i}{\rm ad}^{\bf c}(z_0g_{\ast}^{-1}v_x)}\circ(g^{\bf c}_{\ast})^{-1}
\right),}
\end{array}$$
where $g_{\ast}^{\bf c}$ (resp. ${\rm ad}^{\bf c}$) is the complexification 
of $g_{\ast}$ (resp. ${\rm ad}$).  
Also, we call 
${\rm Ker}(D^{co}_{z_0v_x}-z_0D^{si}_{z_0v_x}\circ A_x^{{\bf c}})$ 
the {\it foccal space} of the complex focal radius $z_0$ and its complex 
dimension the {\it multiplicity} of the complex focal radius $z_0$, 
In [Koi3], it was shown that, in the case where $M$ is of class $C^{\omega}$, 
complex focal radii of $M$ at $x$ indicate the positions of focal points 
of the extrinsic complexification $M^{\bf c}(\hookrightarrow 
G^{\bf c}/K^{\bf c})$ of $M$ along the complexified geodesic 
$\gamma_{v_x}^{\bf c}$, where $G^{\bf c}/K^{\bf c}$ is the anti-Kaehlerian 
symmetric space associated with $G/K$.  
See [Koi3] (also [Koi10]) about the detail of the definition of 
the extrinsic complexification.  
Denote by ${\cal CFR}_x$ the set of all complex focal radii of $M$ at $x$.  
If ${\cal CFR}_x$ is independent of the choice of $x$, then $M$ is called 
a {\it complex equifocal hypersurface}.  
Here we note that we should call such a hypersurface an equi-complex focal 
hypersurface but, for simplicity, we call it a complex equifocal 
hypersurface.  
This notion is the hypersurface version of a complex equifocal submanifold 
defined in [Koi2].  

Next we recall the notion of an anti-Kaehlerian equifocal hypersurface in 
an anti-Kaehlerian symmetric space.  
Let $J$ be a parallel complex structure on an even dimensional 
pseudo-Riemannian manifold $(M,\langle\,\,,\,\,\rangle)$ of half index.  If 
$\langle JX,JY\rangle=-\langle X,Y\rangle$ holds for every $X,\,Y\in TM$, 
then $(M,\langle\,\,,\,\,\rangle,J)$ is called an 
{\it anti-Kaehlerian manifold}.  
Let $N=G/K$ be a symmetric space of non-compact type and 
$G^{\bf c}/K^{\bf c}$ the anti-Kaehlerian symmetric space associated with 
$G/K$.  See [Koi3] about the anti-Kaehlerian structure of 
$G^{\bf c}/K^{\bf c}$.  
Let $f$ be an isometric immersion of 
an anti-Kaehlerian manifold $(M,\langle\,\,,\,\,\rangle,J)$ into 
$G^{\bf c}/K^{\bf c}$.  
If $\widetilde J\circ f_{\ast}=f_{\ast}\circ J$, then 
$M$ is called an {\it anti-Kaehlerian submanifold} immersed by $f$.  
Let $A$ be the shape tensor of $M$.  
We have $A_{\widetilde Jv}X=A_v(JX)=J(A_vX)$,
where $X\in TM$ and $v\in T^{\perp}M$.  
If $A_vX=aX+bJX$ ($a,b\in{\bf R}$), then $X$ is called a 
$J$-{\it eigenvector for} $a+b{\bf i}$.  Let $\{e_i\}_{i=1}^n$ be an 
orthonormal system of $T_xM$ such that $\{e_i\}_{i=1}^n\cup\{Je_i\}_{i=1}^n$ 
is an orthonormal base of $T_xM$.  We call such an orthonormal system 
$\{e_i\}_{i=1}^n$ a $J$-{\it orthonormal base} of $T_xM$.  
If there exists a $J$-orthonormal base consisting of $J$-eigenvectors of 
$A_v$, then we say that $A_v$ {\it is diagonalizable with respect to an} 
$J$-{\it orthonormal base}.  
Then we set ${\rm Tr}_JA_v:=\sum\limits_{i=1}^n\lambda_i$ as 
$A_ve_i=({\rm Re}\,\lambda_i)e_i+({\rm Im}\,\lambda_i)Je_i\,\,\,
(i=1,\cdots,n)$.  We call this quantity the $J$-{trace of} $A_v$.  
If, for each unit normal vector $v\in M$, the shape operator $A_v$ is 
diagonalizable with respect to a $J$-orthonormal tangent base, 
if the normal Jacobi operator $R(v)$ preserves the tangent space $T_xM$ 
($x\,:\,$the base point of $v$) invariantly and if $A_v$ and $R(v)$ 
commute, then we call $M$ a {\it curvature-adapted anti-Kaehlerian 
submanifold}, where $R$ is the curvature tensor of $G^{\bf c}/K^{\bf c}$.  
Assume that $M$ is an anti-Kaehlerian hypersurface (i.e., 
${\rm codim}\,M=2$) and that it is orientable.  
Denote by $\exp^{\perp}$ the normal exponential map of $M$.  
Fix a global parallel orthonormal normal base $\{v,Jv\}$ of $M$.  
If $\exp^{\perp}(av_x+bJv_x)$ is a focal point of $(M,x)$, then we call 
the complex number $a+b{\bf i}$ a {\it complex focal radius along the 
geodesic} $\gamma_{v_x}$.  
Assume that the number (which may be $0$ and $\infty$) of distinct complex 
focal radii along the geodesic $\gamma_{v_x}$ 
is independent of the choice of $x\in M$.  
Furthermore assume that the number is not equal to $0$.  
Let $\{r_{i,x}\,\vert\,i=1,2,\cdots\}$ be the set of 
all complex focal radii along $\gamma_{v_x}$, where 
$\vert r_{i,x}\vert\,<\,\vert r_{i+1,x}\vert$ or 
"$\vert r_{i,x}\vert=\vert r_{i+1,x}\vert\,\,\&\,\,{\rm Re}\,r_{i,x}
>{\rm Re}\,r_{i+1,x}$" or 
"$\vert r_{i,x}\vert=\vert r_{i+1,x}\vert\,\,\&\,\,
{\rm Re}\,r_{i,x}={\rm Re}\,r_{i+1,x}\,\,\&\,\,
{\rm Im}\,r_{i,x}=-{\rm Im}\,r_{i+1,x}<0$".  
Let $r_i$ ($i=1,2,\cdots$) be complex-valued functions on $M$ defined by 
assigning $r_{i,x}$ to each $x\in M$.  We call 
this function $r_i$ the $i$-{\it th complex focal radius function 
for} $\widetilde v$.  
If the number of distinct complex focal radii along 
$\gamma_{v_x}$ is independent of the choice of $x\in M$, 
complex focal radius functions for $v$ 
are constant on $M$ and they have constant multiplicity, 
then $M$ is called an {\it anti-Kaehlerian equifocal hypersurface}.  
We ([Koi3]) showed the following fact.  

\vspace{0.2truecm}

\noindent
{\bf Fact 3.} {\sl Let $M$ be a complete (embedded) $C^{\omega}$-hypersurface 
in $G/K$.  Then $M$ is complex equifocal if and only if $M^{\bf c}$ is 
anti-Kaehler equifocal.}

\vspace{0.2truecm}

Next we recall the notion of an anti-Kaehlerian isoparametric hypersurface 
in an infinite dimensional anti-Kaehlerian space.  
Let $f$ be an isometric immersion of an anti-Kaehlerian Hilbert manifold 
$(M,\langle\,\,,\,\,\rangle,J)$ into an infinite dimensional anti-Kaehlerian 
space $(V,\langle\,\,,\,\,\rangle,\widetilde J)$.  See Section 5 of [Koi3] 
about the definitions of an anti-Kaehlerian Hilbert manifold and 
an infinite dimensional anti-Kaehlerian space.  
If $\widetilde J\circ f_{\ast}=f_{\ast}\circ J$ holds, then we call $M$ 
an {\it anti-Kaehlerian Hilbert submanifold in} 
$(V,\langle\,\,,\,\,\rangle,\widetilde J)$ {\it immersed by} $f$.  If $M$ is 
of finite codimension and there exists 
an orthogonal time-space decomposition $V=V_-\oplus V_+$ such that 
$\widetilde JV_{\pm}=V_{\mp},\,\,(V,\langle\,\,,\,\,\rangle_{V_{\pm}})$ is 
a Hilbert space, the distance topology associated with 
$\langle\,\,,\,\,\rangle_{V_{\pm}}$ 
coincides with the original topology of $V$ and, for each 
$v\in T^{\perp}M$, the shape operator $A_v$ is a compact operator 
with respect to $f^{\ast}\langle\,\,,\,\,\rangle_{V_{\pm}}$, then we call $M$ 
a {\it anti-Kaehlerian Fredholm submanifold} (rather than 
{\it anti-Kaehlerian Fredholm Hilbert submanifold}).  Let $(M,\langle\,\,,\,\,
\rangle,J)$ be an orientable anti-Kaehlerian Fredholm hypersurface in 
an anti-Kaehlerian space $(V,\langle\,\,,\,\,\rangle,\widetilde J)$ and 
$A$ be the shape tensor of $(M,\langle\,\,,\,\,\rangle,J)$.  
Fix a global unit normal vector field $v$ of $M$.  
If there exists $X(\not=0)\in T_xM$ with 
$A_{v_x}X=aX+bJX$, then we call the complex number $a+b{\bf i}$ a 
$J$-{\it eigenvalue of} $A_{v_x}$ (or a {\it complex principal curvature of} 
$M$ {\it at} $x$) and call $X$ a $J$-{\it eigenvector of} $A_{v_x}$ {\it for} 
$a+b{\bf i}$.  
Here we note that this relation is rewritten as $A^{\bf c}_{v_x}
X^{(1,0)}=(a+b{\bf i})X^{(1,0)}$, where $X^{(1,0)}:=\frac12(X-{\bf i}JX)$.  
Also, we call the space of all $J$-eigenvectors of $A_{v_x}$ for 
$a+b\sqrt{-1}$ a $J$-{\it eigenspace of} $A_{v_x}$ {\it for} $a+b{\bf i}$.  
We call the set of all $J$-eigenvalues of $A_{v_x}$ the $J$-{\it spectrum of} 
$A_{v_x}$ and denote it by ${\rm Spec}_JA_{v_x}$.  
${\rm Spec}_JA_{v_x}\setminus\{0\}$ is described as follows:
$${\rm Spec}_JA_{v_x}\setminus\{0\}
=\{\lambda_i\,\vert\,i=1,2,\cdots\}$$
$$\left(
\begin{array}{c}
\displaystyle{\vert\lambda_i\vert>\vert\lambda_{i+1}\vert\,\,\,{{\rm or}}
\,\,\,{\rm "}\vert\lambda_i\vert=\vert\lambda_{i+1}\vert\,\,\&\,\,
{{\rm Re}}\,\lambda_i>{{\rm Re}}\,\lambda_{i+1}{\rm "}}\\
\displaystyle{{{\rm or}}\,\,\,"\vert\lambda_i\vert=\vert\lambda_{i+1}\vert 
\,\,\&\,\,
{{\rm Re}}\,\lambda_i={{\rm Re}}\,\lambda_{i+1}\,\,\&\,\,
{{\rm Im}}\,\lambda_i=-{{\rm Im}}\,\lambda_{i+1}>0"}
\end{array}
\right).$$
Also, the $J$-eigenspace for each $J$-eigenvalue of $A_{v_x}$ other than $0$ 
is of finite dimension.  
We call the $J$-eigenvalue $\lambda_i$ 
the $i$-{\it th complex principal curvature of} $M$ {\it at} $x$.  
Assume that the number 
(which may be $\infty$) of distinct complex principal curvatures of $M$
is constant over $M$.  
Then we can define functions $\widetilde{\lambda}_i$ ($i=1,2,\cdots$) on $M$ 
by assigning the $i$-th complex principal curvature of $M$ at $x$ 
to each $x\in M$.  We call this function 
$\widetilde{\lambda}_i$ the $i$-{\it th complex principal curvature function} 
of $M$.  If the number of distinct complex principal curvatures of $M$ 
is constant over $M$, each complex principal curvature 
function is constant over $M$ and it has constant multiplicity, then 
we call $M$ an {\it anti-Kaehler isoparametric hypersurface}.  
Let $\{e_i\}_{i=1}^{\infty}$ be an orthonormal system of 
$(T_xM,\,\,\langle\,\,,\,\,\rangle_x)$.  If 
$\{e_i\}_{i=1}^{\infty}\cup\{Je_i\}_{i=1}^{\infty}$ is an orthonormal base 
of $T_xM$, then we call $\{e_i\}_{i=1}^{\infty}$ a $J$-{\it orthonormal base}. 
If there exists a $J$-orthonormal base consisting of $J$-eigenvectors of 
$A_{v_x}$, then $A_{v_x}$ is said to {\it be diagonalized with respect to the} 
$J$-{\it orthonormal base}.  If $M$ is anti-Kaehlerian isoparametric and, 
for each $x\in M$, the shape operator $A_{v_x}$ is diagonalized with respect 
to an $J$-orthonormal base, then we call $M$ a {\it proper anti-Kaehlerian 
isoparametric hypersurface}.  

In [Koi2], we defined the notion of the parallel transport map for 
a semi-simple Lie group $G$ as a pseudo-Riemannian submersion of a pseudo-Hilbert space 
$H^0([0,1],\mathfrak g)$ onto $G$.  See [Koi2] in detail.  
Also, in [Koi3], we defined the notion of the parallel transport map for the 
complexification $G^{\bf c}$ of a semi-simple Lie group $G$ as 
an anti-Kaehlerian submersion of an infinite dimensional anti-Kaehlerian 
space $H^0([0,1],\mathfrak g^{\bf c})$ onto $G^{\bf c}$.  See [Koi3] 
in detail.  Let $G/K$ be a symmetric space of non-compact type and 
$\phi:H^0([0,1],\mathfrak g^{\bf c})\to G^{\bf c}$ the parallel transport map 
for $G^{\bf c}$ and $\pi:G^{\bf c}\to G^{\bf c}/K^{\bf c}$ the natural 
projection.  We [Koi3] showed the following fact.  

\vspace{0.2truecm}

\noindent
{\bf Fact. 4.} {\sl Let $M$ be a complete anti-Kaehlerian hypersurface in 
an anti-Kaehlerian symmetric space $G^{\bf c}/K^{\bf c}$.  Then $M$ is 
anti-Kaehlerian equifocal if and only if each component of 
$(\pi\circ\phi)^{-1}(M)$ is anti-Kaehlerian isoparametric.}

\vspace{0.2truecm}

Next we recall the notion of a focal point of non-Euclidean type on the ideal 
boundary $N(\infty)$ of a hypersurface $M$ in a Hadamard manifold $N$ 
which was introduced in [Koi7] for a submanifold of general codimension.  
Assume that $M$ is orientable.  Let $v$ be a unit normal vector field of $M$ and 
$\gamma_{v_x}:[0,\infty)\to N$ the normal geodesic of $M$ of direction $v_x$.  
If there exists a $M$-Jacobi field 
$Y$ along $\gamma_{v_x}$ satisfying 
$\lim\limits_{t\to\infty}\frac{\vert\vert Y(t)\vert\vert}{t}=0$, then we call 
$\gamma_{v_x}(\infty)\,(\in N(\infty))$ 
a {\it focal point} of $M$ {\it on the ideal boundary} $N(\infty)$ 
{\it along} $\gamma_{v_x}$, 
where $\gamma_{v_x}(\infty)$ is the asymptotic class of $\gamma_{v_x}$.  
Also, if there exists a $M$-Jacobi field $Y$ along $\gamma_{v_x}$ satisfying 
$\lim\limits_{t\to\infty}\frac{\vert\vert Y(t)\vert\vert}{t}=0$ and 
${\rm Sec}(v_x,Y(0))\not=0$, then we call $\gamma_{v_x}(\infty)$ a 
{\it focal point of non-Euclidean type of} $M$ {\it on} $N(\infty)$ 
{\it along} $\gamma_{v_x}$, where ${\rm Sec}(v_x,Y(0))$ is the 
sectional curvature for the $2$-plane spanned by $v_x$ and $Y(0)$.  
If, for any point $x$ of $M$, 
$\gamma_{v_x}(\infty)$ and $\gamma_{-v_x}(\infty)$ are not a focal point of 
non-Euclidean type of $M$ on $N(\infty)$, then we say 
that $M$ {\it has no focal point of non-Euclidean type on the ideal boundary} 
$N(\infty)$.  
According to Theorem 1 of [Koi3] and Theorem A of [Koi7], 
we have the following fact.  

\vspace{0.2truecm}

\noindent
{\bf Fact 5.} {\sl Let $M$ be a curvature-adapted and 
isoparametric $C^{\omega}$-hypersurface in a symmetric space $N:=G/K$ of 
non-compact type.  
Then the following conditions ${\rm (i)}$ and ${\rm (ii)}$ are equivalent:

{\rm (i)} $M$ has no focal point of non-Euclidean type on the ideal boundary 
$N(\infty)$.

{\rm(ii)} each component of $(\pi\circ\phi)^{-1}(M^{\bf c})$ is proper 
anti-Kaehlerian isoparametric.}


\section{Proof of Theorem A}
In this section, we shall prove Theorem A.  
Let $M$ be a curvature-adapted isoparametric hypersurface in a simply 
connected symmetric space $G/K$ of compact type, $v$ a unit normal vector 
field of $M$ and $C(\subset T^{\perp}_xM)$ the Coxeter domain 
(i.e., the fundamental domain (containing $0$) of the Coxeter group of $M$ at 
$x$).  The boundary $\partial C$ of $C$ consists of two points and it is 
described as $\partial C=\{r_1v_x,r_2v_x\}$ ($r_2<0<r_1$).  
We may assume that $\vert r_1\vert\leq\vert r_2\vert$ by replacing $v$ 
with $-v$ if necessary.  
Note that the set ${\cal FR}_M$ of all focal radii of $M$ is equal to 
$\{kr_1+(1-k)r_2\,\vert\,k\in{\Bbb Z}\}$.  
Set $F_i:=\{\gamma_{v_x}(r_i)\,\vert\,x\in M\}$ ($i=1,2$), 
which are all of focal submanifolds of $M$.  
The hypersurface $M$ is the $r_i$-tube over $F_i$ ($i=1,2$).  
Let $\pi$ be the natural projection of $G$ onto $G/K$ and $\phi$ the parallel 
transport map for $G$.  
Let $\widetilde M$ be a component of $(\pi\circ\phi)^{-1}(M)$, which is an 
isoparametric hypersurface in $H^0([0,1],\mathfrak g)$.  
The set ${\cal PC}_{\widetilde M}$ of all principal curvatures other than zero 
of $\widetilde M$ is equal to $\{\frac{1}{kr_1+(1-k)r_2}\,\vert\,
k\in{\Bbb Z}\}$.  Set $\lambda_{2k-1}:=\frac{1}{kr_1+(1-k)r_2}$ 
($k=1,2,\cdots$) and $\lambda_{2k}:=\frac{1}{-(k-1)r_1+kr_2}$ 
($k=1,2,\cdots$).  Then we have 
$\vert\lambda_{i+1}\vert<\vert\lambda_i\vert$ or $\lambda_i=-\lambda_{i+1}>0$ 
for any $i\in{\Bbb N}$.  
Denote by $m_i$ the multiplicity of $\lambda_i$.  
Denote by $A$ (resp. $\widetilde A$) the shape operator of $M$ for $v$ 
(resp. $\widetilde M$ for $v^L$), where $v^L$ is the horizontal lift of $v$ 
to $\widetilde M$ with respect to $\pi\circ\phi$.  
Fix $r_0\in{\cal FR}_M$.  
The focal map $f_{r_0}\,:\,M\to G/K$ is defined by 
$f_{r_0}(x):=\gamma_{v_x}(r_0)$ ($x\in M$).  
Let $F:=f_{r_0}(M)$, which is either $F_1$ or $F_2$.  
Denote by $A^F$ the shape tensor of $F$ and $\psi_t$ the geodesic flow of $G/K$.  

\vspace{0.5truecm}

\centerline{
\unitlength 0.1in
\begin{picture}( 32.7900, 24.7800)( 15.9100,-29.6900)
%
\special{pn 8}%
\special{ar 3052 1970 784 752  0.0000000 6.2831853}%
%
\special{pn 20}%
\special{sh 1}%
\special{ar 3836 1962 10 10 0  6.28318530717959E+0000}%
\special{sh 1}%
\special{ar 3836 1962 10 10 0  6.28318530717959E+0000}%
%
\special{pn 20}%
\special{sh 1}%
\special{ar 3068 1208 10 10 0  6.28318530717959E+0000}%
\special{sh 1}%
\special{ar 3068 1208 10 10 0  6.28318530717959E+0000}%
%
\special{pn 20}%
\special{sh 1}%
\special{ar 3068 2724 10 10 0  6.28318530717959E+0000}%
\special{sh 1}%
\special{ar 3068 2724 10 10 0  6.28318530717959E+0000}%
%
\special{pn 20}%
\special{sh 1}%
\special{ar 2258 1962 10 10 0  6.28318530717959E+0000}%
\special{sh 1}%
\special{ar 2258 1962 10 10 0  6.28318530717959E+0000}%
%
\special{pn 20}%
\special{sh 1}%
\special{ar 3604 1438 10 10 0  6.28318530717959E+0000}%
\special{sh 1}%
\special{ar 3604 1438 10 10 0  6.28318530717959E+0000}%
%
\special{pn 20}%
\special{sh 1}%
\special{ar 2490 1438 10 10 0  6.28318530717959E+0000}%
\special{sh 1}%
\special{ar 2490 1438 10 10 0  6.28318530717959E+0000}%
%
\special{pn 20}%
\special{sh 1}%
\special{ar 2490 2484 10 10 0  6.28318530717959E+0000}%
\special{sh 1}%
\special{ar 2490 2484 10 10 0  6.28318530717959E+0000}%
%
\special{pn 20}%
\special{sh 1}%
\special{ar 3604 2500 10 10 0  6.28318530717959E+0000}%
\special{sh 1}%
\special{ar 3604 2500 10 10 0  6.28318530717959E+0000}%
%
\special{pn 8}%
\special{pa 3692 1962}%
\special{pa 3980 1954}%
\special{fp}%
%
\special{pn 8}%
\special{pa 3516 1508}%
\special{pa 3716 1346}%
\special{fp}%
%
\special{pn 8}%
\special{pa 3068 1416}%
\special{pa 3068 1022}%
\special{fp}%
%
\special{pn 8}%
\special{pa 3060 2578}%
\special{pa 3060 2938}%
\special{fp}%
%
\special{pn 8}%
\special{pa 3516 2416}%
\special{pa 3676 2578}%
\special{fp}%
%
\special{pn 8}%
\special{pa 2578 1516}%
\special{pa 2426 1376}%
\special{fp}%
%
\special{pn 8}%
\special{pa 2474 1970}%
\special{pa 2058 1970}%
\special{fp}%
%
\special{pn 8}%
\special{pa 2562 2424}%
\special{pa 2418 2538}%
\special{fp}%
\put(40.1200,-19.1500){\makebox(0,0)[lt]{$F_1$}}%
\put(30.1100,-29.6900){\makebox(0,0)[lt]{$F_1$}}%
\put(20.4100,-19.2300){\makebox(0,0)[rt]{$F_1$}}%
\put(30.1900,-9.9900){\makebox(0,0)[lb]{$F_1$}}%
\put(37.3200,-12.9900){\makebox(0,0)[lt]{$F_2$}}%
\put(36.8400,-25.3100){\makebox(0,0)[lt]{$F_2$}}%
\put(23.8500,-25.4600){\makebox(0,0)[rt]{$F_2$}}%
\put(24.0100,-13.1500){\makebox(0,0)[rt]{$F_2$}}%
%
\special{pn 13}%
\special{pa 3556 1730}%
\special{pa 3948 1492}%
\special{fp}%
%
\special{pn 13}%
\special{pa 3340 2500}%
\special{pa 3566 2886}%
\special{fp}%
%
\special{pn 13}%
\special{pa 2770 2492}%
\special{pa 2544 2878}%
\special{fp}%
%
\special{pn 13}%
\special{pa 3556 2184}%
\special{pa 3948 2424}%
\special{fp}%
%
\special{pn 13}%
\special{pa 2506 2192}%
\special{pa 2094 2402}%
\special{fp}%
%
\special{pn 13}%
\special{pa 2506 1708}%
\special{pa 2114 1470}%
\special{fp}%
%
\special{pn 13}%
\special{pa 3324 1492}%
\special{pa 3540 1130}%
\special{fp}%
%
\special{pn 13}%
\special{pa 2778 1476}%
\special{pa 2562 1116}%
\special{fp}%
\put(35.8800,-10.6100){\makebox(0,0)[lt]{$M$}}%
\put(39.8000,-14.5300){\makebox(0,0)[lt]{$M$}}%
\put(39.7200,-23.7700){\makebox(0,0)[lt]{$M$}}%
\put(35.5600,-29.1500){\makebox(0,0)[lt]{$M$}}%
\put(24.7300,-29.0800){\makebox(0,0)[lt]{$M$}}%
\put(20.4900,-23.8400){\makebox(0,0)[rt]{$M$}}%
\put(20.5700,-14.0700){\makebox(0,0)[rt]{$M$}}%
\put(25.4600,-9.4600){\makebox(0,0)[rt]{$M$}}%
\put(41.0100,-29.2300){\makebox(0,0)[lt]{$\gamma_{v_x}$}}%
%
\special{pn 20}%
\special{sh 1}%
\special{ar 3740 1608 10 10 0  6.28318530717959E+0000}%
\special{sh 1}%
\special{ar 3740 1608 10 10 0  6.28318530717959E+0000}%
\put(38.2000,-16.1500){\makebox(0,0)[lt]{$x$}}%
%
\special{pn 20}%
\special{sh 1}%
\special{ar 4326 1324 10 10 0  6.28318530717959E+0000}%
\special{sh 1}%
\special{ar 4326 1324 10 10 0  6.28318530717959E+0000}%
%
\special{pn 20}%
\special{pa 4190 1122}%
\special{pa 4486 1576}%
\special{fp}%
\put(43.8900,-13.5300){\makebox(0,0)[lb]{$0$}}%
%
\special{pn 8}%
\special{pa 4566 1100}%
\special{pa 4274 1248}%
\special{fp}%
\special{sh 1}%
\special{pa 4274 1248}%
\special{pa 4342 1236}%
\special{pa 4322 1224}%
\special{pa 4324 1200}%
\special{pa 4274 1248}%
\special{fp}%
\put(45.8100,-10.0700){\makebox(0,0)[lt]{$C$}}%
\put(47.5800,-22.0700){\makebox(0,0)[lt]{$T^{\perp}_xM$}}%
%
\special{pn 8}%
\special{pa 4182 1138}%
\special{pa 3612 1446}%
\special{da 0.070}%
%
\special{pn 8}%
\special{pa 4494 1576}%
\special{pa 3836 1946}%
\special{da 0.070}%
%
\special{pn 20}%
\special{pa 4342 1324}%
\special{pa 4446 1476}%
\special{fp}%
\special{sh 1}%
\special{pa 4446 1476}%
\special{pa 4424 1410}%
\special{pa 4416 1432}%
\special{pa 4392 1432}%
\special{pa 4446 1476}%
\special{fp}%
%
\special{pn 20}%
\special{sh 1}%
\special{ar 4494 1584 10 10 0  6.28318530717959E+0000}%
\special{sh 1}%
\special{ar 4494 1584 10 10 0  6.28318530717959E+0000}%
%
\special{pn 20}%
\special{sh 1}%
\special{ar 4502 1576 10 10 0  6.28318530717959E+0000}%
\special{sh 1}%
\special{ar 4502 1576 10 10 0  6.28318530717959E+0000}%
%
\special{pn 20}%
\special{sh 1}%
\special{ar 4198 1122 10 10 0  6.28318530717959E+0000}%
\special{sh 1}%
\special{ar 4198 1122 10 10 0  6.28318530717959E+0000}%
%
\special{pn 8}%
\special{pa 3836 568}%
\special{pa 4870 2154}%
\special{fp}%
%
\special{pn 20}%
\special{sh 1}%
\special{ar 4798 2046 10 10 0  6.28318530717959E+0000}%
\special{sh 1}%
\special{ar 4798 2046 10 10 0  6.28318530717959E+0000}%
%
\special{pn 20}%
\special{sh 1}%
\special{ar 3876 630 10 10 0  6.28318530717959E+0000}%
\special{sh 1}%
\special{ar 3876 630 10 10 0  6.28318530717959E+0000}%
\put(45.5700,-16.2300){\makebox(0,0)[lb]{$r_1v_x$}}%
\put(42.1300,-11.0700){\makebox(0,0)[lb]{$r_2v_x$}}%
\put(48.3800,-20.7700){\makebox(0,0)[lb]{$(2r_1-r_2)v_x$}}%
\put(44.5300,-13.6100){\makebox(0,0)[lt]{$v_x$}}%
%
\special{pn 8}%
\special{pa 4070 2932}%
\special{pa 3524 2570}%
\special{fp}%
\special{sh 1}%
\special{pa 3524 2570}%
\special{pa 3568 2624}%
\special{pa 3568 2598}%
\special{pa 3590 2590}%
\special{pa 3524 2570}%
\special{fp}%
\put(39.2400,-6.6100){\makebox(0,0)[lb]{$(2r_2-r_1)v_x$}}%
\end{picture}%
\hspace{0truecm}}

\vspace{1truecm}

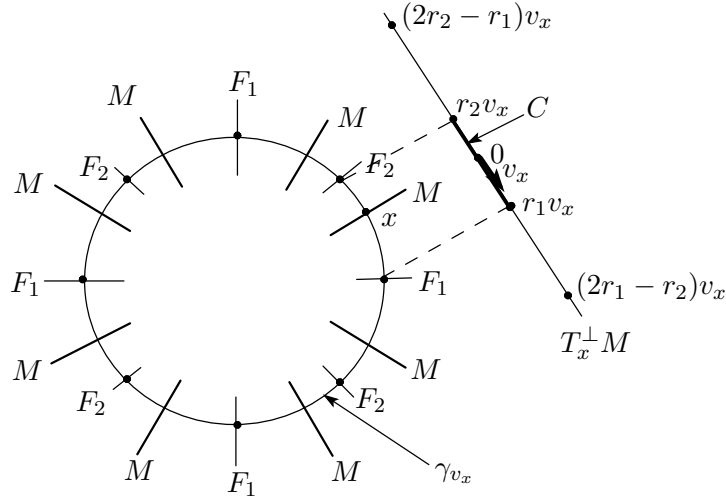
\centerline{{\bf Fig. 2.}}

\vspace{0.5truecm}

\noindent
{\it Proof of Theorem A.} 
Define a set $S_x$ by 
$$S_x:=\{(\lambda,\mu)\in{\rm Spec}A_x\times{\rm Spec}R(v_x)
\,\vert\,{\rm Ker}(A_x-\lambda I)\cap{\rm Ker}(R(v_x)-\mu I)\not=\{0\}\}.$$
Since $M$ is curvature adapted, we have 
$$T_xM=\displaystyle{\mathop{\oplus}_{(\lambda,\mu)\in S_x}
({\rm Ker}(A_x-\lambda I)}\cap{\rm Ker}(R(v_x)-\mu I)).$$
Define a distribution $D$ on $M$ by 
$D_x:=\displaystyle{\mathop{\oplus}_{(\lambda,\mu)\in S_{r_0}^x}(
{\rm Ker}(A_x-\lambda I)\cap{\rm Ker}(R(v_x)-\mu I))}$ and $D^{\perp}$ 
the orthogonal complementary distribution of $D$ in $TM$.  
Let $X\in{\rm Ker}(A_x-\lambda I)\cap{\rm Ker}(R(v_x)-\mu I)$ 
($(\lambda,\mu)\in S_{r_0}^x$) and $Y$ be the Jacobi field along 
$\gamma_{r_0v_x}$ with $Y(0)=X$ and $Y'(0)=-A_{r_0v_x}X\,(=-r_0\lambda X)$.  
This Jacobi field $Y$ is described as 
$$Y(s)=\left(\cos(sr_0\sqrt{\mu})-
\frac{\lambda\sin(sr_0\sqrt{\mu})}
{\sqrt{\mu}}\right)P_{\gamma_{r_0v}\vert_{[0,s]}}(X).$$
Since $Y(1)=f_{r_0\ast}X$, we have 
$$f_{r_0\ast}X=\left(\cos(r_0\sqrt{\mu})
-\frac{\lambda\sin(r_0\sqrt{\mu})}
{\sqrt{\mu}}\right)P_{\gamma_{r_0v_x}}(X),\leqno{(3.1)}$$
which is not equal to $0$ because $(\lambda,\mu)\in S_{r_0}^x$.  From this 
relation, we have $T_{f_{r_0}(x)}F=P_{\gamma_{r_0v_x}}(D)$.  
On the other hand, we have 
$$\begin{array}{l}
\displaystyle{
\hspace{0.3truecm}\widetilde{\nabla}_{f_{r_0\ast}X}\psi_{r_0}(v_x)
=\frac{1}{r_0}Y'(1)}\\
\displaystyle{=-\left(\sqrt{\mu}\sin(r_0\sqrt{\mu})
+\lambda\cos(r_0\sqrt{\mu})\right)P_{\gamma_{r_0v_x}}(X).}
\end{array}
\leqno{(3.2)}
$$
From $(3.1)$ and $(3.2)$, we have 
$$A^F_{\psi_{r_0}(v_x)}f_{r_0\ast}X
=-\frac{\mu+\lambda \tau_{r_0}(\mu)}{\lambda-\tau_{r_0}(\mu)}f_{r_0\ast}X.$$
Hence we can derive the following relation:
$${\rm Tr}\,A^F_{\psi_{r_0}(v_x)}
=-\sum_{(\lambda,\mu)\in S_{r_0}^x}\frac{\mu+\lambda\tau_{r_0}(\mu)}
{\lambda-\tau_{r_0}(\mu)}\times m_{\lambda,\mu},\leqno{(3.3)}$$
where $S_{r_0}^x$ and $m_{\lambda,\mu}$ are as in the statement of 
Theorem A.  
On the other hand, it is not difficult to show the existence of a transnormal function on 
$G/K$ having $M$ and $F$ as a regular level and a singular level, respectively.  
Hence, according to Theorem 1.3 of [Mi], $F$ is austere and hence minimal.  
Therefore, we obtain the desired identity from $(3.3)$.  
\begin{flushright}q.e.d.\end{flushright}

\section{The mean curvature of a proper anti-Kaehlerian Fredholm submanifold}
In this section, we define the notion of a proper anti-Kaehlerian Fredholm 
submanifold and its mean curvature vector.  
Let $M$ be an anti-Kaehlerian Fredholm submanifold in an infinite 
dimensional anti-Kaehlerian space $V$ and $A$ be the shape tensor of $M$.  
Denote by the same symbol $J$ the complex structures of $M$ and $V$.  
If $A_v$ is diagonalized with respect to 
a $J$-orthonormal base for each unit normal vector $v$ of $M$, then we 
call $M$ a {\it proper anti-Kaehlerian Fredholm submanifold}.  Assume that 
$M$ is such a submanifold.  Let $v$ be a unit normal vector of $M$.  
If the series $\sum\limits_{i=1}^{\infty}m_i\lambda_i$ exists, then we 
call it the $J$-{\it trace of} $A_v$ and denote it by ${\rm Tr}_JA_v$, 
where $\{\lambda_i\,\vert\,i=1,2,\cdots\}={\rm Spec}_JA_v\setminus\{0\}$ 
($\lambda_i$'s are ordered as stated in Section 2) and 
$m_i=\frac12{\rm dim}{\rm Ker}(A_v-\lambda_iI)$ ($i=1,2,\cdots$), where 
$\lambda_iI$ means $({\rm Re}\,\lambda_i)I+({\rm Im}\,\lambda_i)J$.  
Note that, if $\sharp({\rm Spec}_JA_v)$ is 
finite, then we promise $\lambda_i=0$ and $m_i=0$ 
($i>\sharp({\rm Spec}_JA_v\setminus\{0\})$), where $\sharp(\cdot)$ is 
the cardinal number of ($\cdot$).  
Define a normal vector field $H$ of $M$ by $\langle H_x,v\rangle
={\rm Tr}_JA_v$ ($x\in M,\,\,v\in T_x^{\perp}M$).  We call $H$ the 
{\it mean curvature vector of} $M$.  

Let $G/K$ be a symmetric space of non-compact type and 
$\phi:H^0([0,1],\mathfrak g^{\bf c})\to G^{\bf c}$ be 
the parallel transport map for the complexification $G^{\bf c}$ of $G$ and 
$\pi$ be the natural projection of $G^{\bf c}$ onto the anti-Kaehlerian 
symmetric space $G^{\bf c}/K^{\bf c}$.  We have the following fact, 
which will be used in the proof of Theorem B in the next section.  

\vspace{0.5truecm}

\noindent
{\bf Lemma 4.1.}  {\sl Let $M$ be a curvature-adapted anti-Kaehlerian 
submanifold in $G^{\bf c}/K^{\bf c}$ and $A$ {\rm(}resp. $\widetilde A${\rm)} 
be the shape tensor of $M$ {\rm(}resp. $(\pi\circ\phi)^{-1}(M)${\rm)}.  
Assume that, for each unit normal vector $v$ of $M$ and each $J$-eigenvalue 
$\mu$ of $R(v)$, ${\rm Ker}(A_v-\sqrt{-\mu}I)\,\cap\,{\rm Ker}
(R(v)-\mu I)=\{0\}$ holds.  
Then the following statements {\rm(i)} and {\rm (ii)} hold:

{\rm(i)} $(\pi\circ\phi)^{-1}(M)$ is a proper anti-Kaehlerian Fredholm 
submanifold.

{\rm (ii)} For each unit normal vector $v$ of $M$, ${\rm Tr}_J
\widetilde A_{v^L}={\rm Tr}_JA_v$ holds, 
where $v^L$ is the horizontal lift of $v$ to $(\pi\circ\phi)^{-1}(M)$ 
and ${\rm Tr}_JA_v$ is the $J$-trace of $A_v$.  
}

\vspace{0.5truecm}

\noindent
{\it Proof.} We can show the statement (i) in terms of Lemmas 9, 12 and 13 
in [Koi3].  By imitating the proof of Theorem C in [Koi2], we can show 
the statement (ii), where we also use the above lemmas in [Koi3].  
\hspace{8.4truecm}q.e.d.

\section{Proofs of Theorems B and C}
In this section, we first prove Theorem B.  
Let $M$ be a curvature-adapted isoparametric $C^{\omega}$-hypersurface 
in a symmetric space $G/K$ of non-compact type.  
Assume that $M$ admits no focal point of non-Euclidean type on the ideal 
boundary of $G/K$.  
Denote by $A$ the shape tensor of $M$ and $R$ the curvature tensor of $G/K$.  
Let $v$ be a unit normal vector field of $M$, which is uniquely extended to 
a unit normal vector field of the extrinsic complexification 
$M^{\bf c}(\subset G^{\bf c}/K^{\bf c})$ of $M$.  
Since $M$ is a curvature-adapted isoparametric hypersurface admitting no 
focal point of non-Euclidean type on the ideal boundary $N(\infty)$, 
it admits a complex focal radius.  
Let $r_0$ be one of complex focal radii of $M$.  
The focal map $f_{r_0}\,:\,M^{\bf c}\to G^{\bf c}/K^{\bf c}$ for $r_0$ is 
defined by $f_{r_0}(x):=\exp^{\perp}(r_0v_x)(=\gamma_{v_x}^{\bf c}(r_0))$ 
($x\in M^{\bf c}$), where $r_0v_x$ means $({\rm Re}r_0)v_x+({\rm Im}r_0)Jv_x$ 
($J\,:\,$ the complex structure of $G^{\bf c}/K^{\bf c}$).  Let 
$F:=f_{r_0}(M^{\bf c})$, which is an anti-Kaehlerian submanifold in 
$G^{\bf c}/K^{\bf c}$ (see Fig. 1).  
Without loss of generality, we may assume $o:=eK\in M$.  
Denote by $\widehat A$ and $A^F$ the shape tensor of $M^{\bf c}$ and $F$, 
respectively.  Let $\psi_t$ be the geodesic flow of $G^{\bf c}/K^{\bf c}$.  
Then we have the following fact.  

\vspace{0.5truecm}

\noindent
{\bf Lemma 5.1.} {\sl For any $x\in M\,(\subset M^{\bf c})$, 
the following relation holds:
$${\rm Tr}_JA^F_{\psi_{\vert r_0\vert}(\frac{r_0}{\vert r_0\vert}v_x)}
=-\frac{r_0}{\vert r_0\vert}\sum_{(\lambda,\mu)\in S_{r_0}^x}
\frac{\mu+\lambda\hat{\tau}_{r_0}(\mu)}
{\lambda-\hat{\tau}_{r_0}(\mu)}
\times m_{\lambda,\mu},$$
where $S_{r_0}^x$ and $m_{\lambda,\mu}$ are as in the statement of 
Theorem B.}

\vspace{0.5truecm}

\noindent
{\it Proof.} Let $S_x:=\{(\lambda,\mu)\in{\rm Spec}A_{v_x}\times{\rm Spec}
R(v_x)\,\vert\,{\rm Ker}(A_{v_x}-\lambda I)\cap{\rm Ker}(R(v_x)-\mu I)
\not=\{0\}\}$.  Since $M$ is curvature adapted, we have 
$T_xM=\displaystyle{\mathop{\oplus}_{(\lambda,\mu)\in S_x}
\left({\rm Ker}(A_x-\lambda I)\cap{\rm Ker}(R(v_x)-\mu I)\right)}$.  
Set $D_x:=\displaystyle{\mathop{\oplus}_{(\lambda,\mu)\in S_{r_0}^x}\left(
{\rm Ker}(A_x-\lambda I)\cap{\rm Ker}(R(v_x)-\mu I)\right)}$ and 
$D^{\perp}_x$ the orthogonal complement of $D_x$ in $T_xM$.  The tangent space 
$T_x(M^{\bf c})$ is identified with the complexification $(T_xM)^{\bf c}$.  
Under this identification, the shape operator $\widehat A_{v_x}$ 
is identified with the complexification $A_x^{\bf c}$ of $A_x$.  
Let $X\in{\rm Ker}(A_x-\lambda I)^{\bf c}\cap
{\rm Ker}(R(v_x)-\mu I)^{\bf c}$ ($(\lambda,\mu)\in S_{r_0}^x$) and $Y$ be the 
Jacobi field along $\gamma_{r_0v_x}$ with $Y(0)=X$ and $Y'(0)=-\hat A_{r_0v_x}
X\,(=-r_0\lambda X=-\lambda\left(({\rm Re}r_0)X+({\rm Im}r_0)JX\right)$), 
where $\gamma_{r_0v_x}$ is the geodesic in $G^{\bf c}/K^{\bf c}$ with 
$\dot{\gamma}_{r_0v_x}(0)=r_0v_x(=({\rm Re}r_0)v_x+({\rm Im}r_0)Jv_x)$.  
This Jacobi field $Y$ is described as 
$$Y(s)=\left(\cos({\bf i}sr_0\sqrt{-\mu})-
\frac{\lambda\sin({\bf i}sr_0\sqrt{-\mu})}
{{\bf i}\sqrt{-\mu}}\right)P_{\gamma_{r_0v_x}\vert_{[0,s]}}(X).$$
Since $Y(1)=f_{r_0\ast}X$, we have 
$$f_{r_0\ast}X=\left(\cos({\bf i}r_0\sqrt{-\mu})
-\frac{\lambda\sin({\bf i}r_0\sqrt{-\mu})}
{{\bf i}\sqrt{-\mu}}\right)P_{\gamma_{r_0v_x}}(X)\leqno{(5.1)}$$
which is not equal to $0$ because $(\lambda,\mu)\in S_{r_0}^x$.  This 
relation implies that $T_{f_{r_0}(x)}F=P_{\gamma_{r_0v_x}}(D_x^{\bf c})$.  
On the other hand, we have 
$$\begin{array}{l}
\displaystyle{
\hspace{0.3truecm}\widetilde{\nabla}_{f_{r_0\ast}X}\psi_{\vert r_0\vert}
(\frac{r_0}{\vert r_0\vert}v_x)=\frac{1}{\vert r_0\vert}Y'(1)}\\
\displaystyle{=-\frac{r_0}{\vert r_0\vert}\left({\bf i}\sqrt{-\mu}
\sin({\bf i}r_0\sqrt{-\mu})+\lambda\cos({\bf i}r_0\sqrt{-\mu})\right)
P_{\gamma_{r_0v_x}}(X).}
\end{array}
\leqno{(5.2)}
$$
From $(5.1)$ and $(5.2)$, we have 
$$A^F_{\psi_{\vert r_0\vert}(\frac{r_0}{\vert r_0\vert}v_x)}f_{r_0\ast}X
=\frac{-\frac{r_0}{\vert r_0\vert}\left(\mu
+\lambda\hat{\tau}_{r_0}(\mu)\right)}
{\lambda-\hat{\tau}_{r_0}(\mu)}
f_{r_0\ast}X.
\leqno{(5.3)}$$
The desired relation follows from this relation.\hspace{6truecm}q.e.d.

\vspace{0.5truecm}

Set $\kappa(\lambda,\mu):=\frac{-\frac{r_0}{\vert r_0\vert}
(\mu+\lambda\hat{\tau}_{r_0}(\mu))}{\lambda-\hat{\tau}_{r_0}(\mu)}$ 
($(\lambda,\mu)\in S_{r_0}^x$).  
Next we prepare the following lemma.  

\vspace{0.5truecm}

\noindent
{\bf Lemma 5.2.} {\sl Let $(\lambda_1,\mu_1)\in S_{r_0}^x$.  Then we have 

{\rm (i)} $\displaystyle{(\exp_{G^{\bf c}}r_0v_x)_{\ast}^{-1}\psi_{\vert r_0
\vert}(\frac{r_0}{\vert r_0\vert}v_x)=\frac{r_0}{\vert r_0\vert}v_x}$, 
where $\exp_{G^{\bf c}}$ is the exponential map of $G^{\bf c}$,

{\rm (ii)} $\displaystyle{(\exp_{G^{\bf c}}r_0v_x)_{\ast}^{-1}\left(
{\rm Ker}(A^F_{\psi_{\vert r_0\vert}(\frac{r_0}{\vert r_0\vert}v_x)}-
\kappa(\lambda_1,\mu_1)I)\right)}$

\hspace{1.2truecm} $\displaystyle{=\mathop{\oplus}_{(\lambda,\mu)\in S_{r_0}^x
(\lambda_1,\mu_1)}
\left({\rm Ker}(A_{v_x}-\lambda I)^{\bf c}\cap{\rm Ker}(R(v_x)-\mu I)^{\bf c}
\right)}$, \newline
where $S_{r_0}^x(\lambda_1,\mu_1)=\{(\lambda,\mu)\in S_{r_0}^x\,\vert\,
\kappa(\lambda,\mu)=\kappa(\lambda_1,\mu_1)\}$, 

{\rm (iii)} if $\lambda_1\not=\pm\sqrt{-\mu_1}$, then 
$\kappa(\lambda_1,\mu_1)\not=\pm\frac{r_0}{\vert r_0\vert}\sqrt{-\mu_1}$.}

\vspace{0.5truecm}

\noindent
{\it Proof.}  The relation of (i) is trivial.  Let $(\lambda,\mu)\in 
S_{r_0}^x(\lambda_1,\mu_1)$.  The restriction \newline
$f_{r_0\ast}\vert_{{\rm Ker}(A_{v_x}-\lambda I)^{\bf c}\cap
{\rm Ker}(R(v_x)-\mu I)^{\bf c}}$ 
of $f_{r_0\ast}$ is equal to $P_{\gamma_{r_0v_x}}\vert_
{{\rm Ker}(A_{v_x}-\lambda I)^{\bf c}\cap{\rm Ker}(R(v_x)-\mu I)^{\bf c}}$ up 
to constant multiple by $(5.1)$.  Also, we have $P_{\gamma_{r_0v_x}}=
(\exp_{G^{\bf c}}r_0v_x)_{\ast}$.  These facts together with $(5.3)$ deduce 
$$\begin{array}{l}
\displaystyle{(\exp_{G^{\bf c}}r_0v_x)_{\ast}
\left({\rm Ker}(A_{v_x}-\lambda I)^{\bf c}\cap{\rm Ker}(R(v_x)-\mu I)^{\bf c}
\right)}\\
\displaystyle{=f_{r_0\ast}\left(
{\rm Ker}(A_{v_x}-\lambda I)^{\bf c}\cap{\rm Ker}(R(v_x)-\mu I)^{\bf c}
\right)}\\
\displaystyle{\subset{\rm Ker}\left(A^F_{\psi_{\vert r_0\vert}(\frac{r_0}
{\vert r_0\vert}v_x)}-\kappa(\lambda_1,\mu_1)I\right).}
\end{array}$$  
From this fact, the relation of (ii) follows.  
Now we shall show the statement (iii).  Let $r_0=a_0+b_0\sqrt{-1}$ 
($a_0,b_0\in{\bf R}$).  Suppose that $\kappa(\lambda_1,\mu_1)
=\pm\frac{r_0}{\vert r_0\vert}\sqrt{-\mu_1}$.  
By squaring both sides of this relation, we have 
$$\left(\hat{\tau}_{r_0}(\mu_1)^2+\mu_1\right)
\left(\lambda_1^2+\mu_1\right)=0.$$
Hence we have $\lambda_1=\pm\sqrt{-\mu_1}$.  
Thus the statement (iii) is shown.  \hspace{2.7truecm} q.e.d.

\vspace{0.5truecm}

Denote by $\hat R$ the curvature tensor of $G^{\bf c}/K^{\bf c}$.  
By using these lemmas, we prove Theorem B.  
According to Lemma 5.1, we have only to show ${\rm Tr}_J
A^F_{\psi_{\vert r_0\vert}(\frac{r_0}{\vert r_0\vert}v_x)}=0$ ($x\in M$).  
In the case where $M$ is homogeneous, we can show this relation by 
imitating the process of the proof of Corollary 1.1 of [HL].  

\vspace{0.5truecm}

\noindent
{\it Simple proof of Theorem B in rank one case.} 
We have only to show 
${\rm Tr}_JA^F_{\psi_{\vert r_0\vert}(\frac{r_0}{\vert r_0\vert}v_x)}=0$.  
Assume that $G/K$ is of rank one.  
Define a complex linear function $\Phi:T^{\perp}_{f_{r_0}(x)}F
\to{\bf C}$ by $\Phi(w)={\rm Tr}_J A^F_w$ ($w\in T^{\perp}_{f_{r_0}(x)}F$).  
Since $M$ is curvature-adapted, we have 
$T_xM=\displaystyle{\mathop{\oplus}_{(\lambda,\mu)\in S_x}
\left({\rm Ker}(A_{v_x}-\lambda I)\cap{\rm Ker}(R(v_x)-\mu I)\right)}$.  
Set 
$$\begin{array}{l}
\displaystyle{
\hat S_{r_0}^y:=\{(\lambda,\mu)\in({\rm Spec}_J\hat A_{v_y})\times
({\rm Spec}_J\hat R(v_y))\,\vert\,{\rm Ker}(\hat A_{v_y}-\lambda I)\cap
{\rm Ker}(\hat R(v_y)-\mu I)\not=\{0\}}\\
\hspace{8truecm}\displaystyle{\,\&\,
\lambda\not=\hat f_{r_0}(\mu)\}}
\end{array}$$
($y\in M^{\bf c}$).  
Define a distribution $\hat D$ on $M^{\bf c}$ by 
$$\hat D_y:=\mathop{\oplus}_{(\lambda,\mu)\in\hat S_{r_0}^y}
\left({\rm Ker}(\hat A_{v_y}-\lambda I)\cap{\rm Ker}(\hat R(v_y)-\mu I)
\right)\quad(y\in M^{\bf c})$$
and ${\hat D}^{\perp}$ the orthogonal complementary 
distribution of $\hat D$ in $T(M^{\bf c})$.  
Also, define a distribution $D$ on $M$ by 
$D_x:=\displaystyle{\mathop{\oplus}_{(\lambda,\mu)\in\hat S_{r_0}^x}
\left({\rm Ker}(A_x-\lambda I)\cap{\rm Ker}(R(v_x)-\mu I)\right)}$ ($x\in M$) 
and $D^{\perp}$ the orthogonal complementary distribution of $D$ in $TM$.  
Under the identification of $T_x(M^{\bf c})$ with $(T_xM)^{\bf c}$, 
$\hat D_x$ is identified with the complexification $(D_x)^{\bf c}$ of $D_x$.  
The focal map $f_{r_0}$ is a submersoin of $M^{\bf c}$ onto $F$ and 
the fibres of $f_{r_0}$ are integral manifolds of ${\hat D}^{\perp}$.  
Let $L$ be the integral manifold of ${\hat D}^{\perp}$ through $x$ 
and set $L_{\bf R}:=L\cap M$.  It is shown that $L$ is the extrinsic 
complexification of $L_{\bf R}$.  
Set 
$Q:=\{\psi_{\vert r_0\vert}(\frac{r_0}{\vert r_0\vert}v_x)\,\vert\,x\in L\}$ 
and $Q_{\bf R}:=\{\psi_{\vert r_0\vert}(\frac{r_0}{\vert r_0\vert}v_x)\,\vert\,
x\in L_{\bf R}\}$.  It is shown that $Q$ is the extrinsic complexification of 
$Q_{\bf R}$ and that $Q$ is a complex hypersurface without geodesic point in 
$T^{\perp}_{f_{r_0}(x)}F$, that is, 
it is not contained in any complex affine hyperplane of 
$T^{\perp}_{f_{r_0}(x)}F$.  According to Lemma 5.1, 
we have 
$$\Phi(\psi_{\vert r_0\vert}(\frac{r_0}{\vert r_0\vert}v_y))=
-\frac{r_0}{\vert r_0\vert}\sum_{(\lambda,\mu)\in S^y_{r_0}}
\frac{\mu+\lambda\hat{\tau}_{r_0}(\mu)}{\lambda-\hat{\tau}_{r_0}(\mu)}\times 
m_{\lambda,\mu}.$$
Let $(\widetilde{\lambda},\widetilde{\mu})$ be a pair of continuous functions 
on $L_{\bf R}$ such that $(\widetilde{\lambda}(y),\widetilde{\mu}(y))\in 
S^y_{r_0}$ for any $y\in L$.  
Since $G/K$ is of rank one, $\widetilde{\mu}$ is constant on $L_{\bf R}$.  
The complex focal radius having 
${\rm Ker}(A_y-\widetilde{\lambda}(y)\,I)\cap{\rm Ker}(R(v_y)-\widetilde{\mu}
(y)\,I)$ as a part of the focal space is the complex number $z_0$ satisfying 
${\rm Ker}(D^{co}_{z_0v_y}-z_0D^{si}_{z_0v_y}\circ A^{\bf c}_y)
\vert_{{\rm Ker}(A_y-\widetilde{\lambda}(y)\,I)\cap
{\rm Ker}(R(v_y)-\widetilde{\mu}(y)\,I)}\not=\{0\}$, that is, 
it is equal to 
$\frac{1}{\sqrt{\widetilde{\mu}(y)}}\arctan\,\frac{\sqrt{\widetilde{\mu}(y)}}
{\widetilde{\lambda}(y)}$, which is independent of the choice of 
$y\in L_{\bf R}$ by the isoparametricness (hence complex equifocality) of 
$M$.  Hence $\widetilde{\lambda}$ is constant on $L_{\bf R}$.  
Therefore $\Phi$ is constant along $Q_{\bf R}$.  
Since $\Phi$ is of class $C^{\omega}$ and $Q_{\bf R}$ is 
a half-dimensional totally real submanifold in $Q$, $\Phi$ is constant 
along $Q$.  
Furthermore, this fact together with the linearity of $\Phi$ imply 
$\Phi\equiv 0$.  In particular, we have 
${\rm Tr}\,A^F_{\psi_{r_0}(v_x)}=0$.
\hspace{2truecm}q.e.d.

\vspace{0.5truecm}

\noindent
{\it Proof of Theorem B (general case).} 
According to Lemma 5.1, we have only to show ${\rm Tr}_J
A^F_{\psi_{\vert r_0\vert}(\frac{r_0}{\vert r_0\vert}v_{x_0})}=0$ 
($x_0\in M$).  
We shall show this relation by investigating 
the focal submanifold of $(\pi\circ\phi)^{-1}(M^{\bf c})$ corresponding to $r_0$, 
where $\phi$ ($:H^0([0,1],\mathfrak g^{\bf c})\to G^{\bf c}$) is the 
parallel transport map for $G^{\bf c}$ and $\pi$ is 
the natural projection of $G^{\bf c}$ onto $G^{\bf c}/K^{\bf c}$.  
Let $\widetilde{M^{\bf c}}$ be the complete extension of 
$(\pi\circ\phi)^{-1}(M^{\bf c})$.  
Let $v^L$ be the horizontal lift of $v$ to $\widetilde M^{\bf c}$.  
Since $\pi\circ\phi$ is an anti-Kaehlerian submersion, the complex focal radii 
of $M^{\bf c}$ (hence $M$) are those of $\widetilde M^{\bf c}$.  
Let $r_0$ be a complex focal radius of $M$ (hence $\widetilde M^{\bf c}$).  
The focal map $\widetilde f_{r_0}$ for $r_0$ is defined by 
$\widetilde f_{r_0}(x)=x+r_0v_x^L$ ($x\in\widetilde M^{\bf c}$).  
Set $\widetilde F:=\widetilde f_{r_0}(\widetilde M^{\bf c})$.  
Denote by $\widetilde A$ (resp. $A^{\widetilde F}$) 
the shape tensor of $\widetilde M^{\bf c}$ (resp. $\widetilde F$).  Let 
${\rm Spec}_J\widetilde A_{v^L_{\hat0}}\setminus\{0\}=\{\lambda_i\,\vert\,
i=1,2,\cdots\}$ ("$\vert\lambda_i\vert>\vert\lambda_{i+1}\vert"\,\,{\rm or}\,\,
"\vert\lambda_i\vert=\vert\lambda_{i+1}\vert\,\,\&\,\,{\rm Re}\lambda_i>
{\rm Re}\lambda_{i+1}"\,\,{\rm or}\,\,"\vert\lambda_i\vert=\vert
\lambda_{i+1}\vert\,\,\&\,\,{\rm Re}\lambda_i={\rm Re}\lambda_{i+1}\,\,\&
\,\,{\rm Im}\lambda_i=-{\rm Im}\lambda_{i+1}>0$").  
The set of all complex focal radii of $M^{\bf c}$ (hence $M$) is equl to 
$\{\frac{1}{\lambda_i}\,\vert\,i=1,2,\cdots\}$.  
We have $r_0=\frac{1}{\lambda_{i_0}}$ for some $i_0$.  
Define a distribution $\widetilde D_i$ ($i=0,1,2,\cdots$) on 
$\widetilde M^{\bf c}$ by $(\widetilde D_0)_u:={\rm Ker}
\widetilde A_{\widetilde v_u^L}$ and $(\widetilde D_i)_u:=
{\rm Ker}(\widetilde A_{\widetilde v_u^L}-\lambda_iI)$ ($i=1,2,\cdots$), where 
$u\in\widetilde M^{\bf c}$.  Since $M$ is 
a curvature-adapted isoparametric submanifold admitting no focal point of 
non-Euclidean type on $N(\infty)$, 
$\widetilde M^{\bf c}$ is proper anti-Kaehlerian isoparametric by Fact 5.  
Therefore, we have $T\widetilde M^{\bf c}=
\displaystyle{\overline{\widetilde D_0\oplus
(\mathop{\oplus}_i\widetilde D_i)}}$ and 
${\rm Spec}_J\widetilde A_{\widetilde v_u^L}$ is independent of the choice of 
$u\in\widetilde M^{\bf c}$.  
Take $u_0\in\widetilde M^{\bf c}$ with $(\pi\circ\phi)(u_0)=x_0$.  
Let $X_i\in(\widetilde D_i)_{u_0}$ 
($i\not=i_0$) and $X_0\in (\widetilde D_0)_{u_0}$.  Then we have 
$\widetilde f_{r_0\ast}X_i=(1-r_0\lambda_i)X_i$ and 
$\widetilde f_{r_0\ast}X_0=X_0$.  Hence we have 
$T_{\widetilde f_{r_0}(u_0)}\widetilde F
=(\widetilde D_0)_{u_0}\oplus(\displaystyle{\mathop{\oplus}_{i\not=i_0}
(\widetilde D_i)_{u_0})}$ and ${\rm Ker}(\widetilde f_{r_0})_{\ast u_0}=
(\widetilde D_{i_0})_{u_0}$, which implies that $\widetilde D_{i_0}$ is 
integrable.  On the other hand, we have 
$\displaystyle{A^{\widetilde F}_{\widetilde{\psi}_{\vert r_0\vert}(\frac{r_0}
{\vert r_0\vert}v^L_{u_0})}\widetilde f_{r_0\ast}X_i=
\frac{\lambda_ir_0}{\vert r_0\vert}X_i}$ and 
$\displaystyle{A^{\widetilde F}_{\widetilde{\psi}_{\vert r_0\vert}(\frac{r_0}
{\vert r_0\vert}v^L_{u_0})}\widetilde f_{r_0\ast}X_0=0}$, 
where $\widetilde{\psi}$ is the geodesic flow of 
$H^0([0,1],\frak g^{\bf c})$.  Therefore, we obtain 
$\displaystyle{A^{\widetilde F}_{\widetilde{\psi}_{\vert r_0\vert}(\frac{r_0}
{\vert r_0\vert}v^L_{u_0})}\widetilde f_{r_0\ast}X_i=
\frac{\lambda_i\vert\lambda_{i_0}\vert}{\lambda_{i_0}-\lambda_i}
\widetilde f_{r_0\ast}X_i}$.  
Hence we have ${\rm Tr}_JA^{\widetilde F}_{\widetilde{\psi}_{\vert r_0\vert}
(\frac{r_0}{\vert r_0\vert}v^L_{u_0})}=
\sum\limits_{i\not=i_0}\frac{\lambda_i\vert\lambda_{i_0}\vert}
{\lambda_{i_0}-\lambda_i}\times m_i$, where $m_i:=\frac12{\rm dim}\,
\widetilde D_i$.  
According to Theorem 2 of [Koi3], each leaf of $\widetilde D_{i_0}$ is 
a complex sphere.  
Let $L$ be the leaf of $\widetilde D_{i_0}$ through $u_0$ and 
$u_0^{\ast}$ be the anti-podal point of $u_0$ in the complex sphere $L$.  
Similarly we can show 
$\displaystyle{{\rm Tr}_JA^{\widetilde F}_{\widetilde{\psi}_{\vert r_0\vert}
(\frac{r_0}{\vert r_0\vert}(\widetilde v^L)_{u_0^{\ast}})}=\sum_{i\not=i_0}
\frac{\lambda_i\vert\lambda_{i_0}\vert}{\lambda_{i_0}-\lambda_i}\times m_i}$.  
Thus we have ${\rm Tr}_JA^{\widetilde F}_{\widetilde{\psi}_{\vert r_0\vert}
(\frac{r_0}{\vert r_0\vert}v^L_{u_0})}=
{\rm Tr}_JA^{\widetilde F}_{\widetilde{\psi}_{\vert r_0\vert}
(\frac{r_0}{\vert r_0\vert}(\widetilde v^L)_{u_0^{\ast}})}$. 
On the other hand, it 
follows from $\widetilde{\psi}_{\vert r_0\vert}
(\frac{r_0}{\vert r_0\vert}(\widetilde v^L)_{u_0^{\ast}})
=-\widetilde{\psi}_{\vert r_0\vert}(\frac{r_0}{\vert r_0\vert}v^L_{u_0})$ 
that ${\rm Tr}_JA^{\widetilde F}_
{\widetilde{\psi}_{\vert r_0\vert}(\frac{r_0}{\vert r_0\vert}v^L_{u_0})}
=-{\rm Tr}_JA^{\widetilde F}_{\widetilde{\psi}_{\vert r_0\vert}
(\frac{r_0}{\vert r_0\vert}(\widetilde v^L)_{u_0^{\ast}})}$.  
Hence we obtain 
$${\rm Tr}_JA^{\widetilde F}_{\widetilde{\psi}_{\vert r_0\vert}
(\frac{r_0}{\vert r_0\vert}v^L_{u_0})}=0. \leqno{(5.4)}$$
It follows from (i) and (ii) of Lemma 5.2 that $F:=f_{r_0}(M^{\bf c})$ is 
a curvature adapted anti-Kaehlerian submanifold.  Also, it follows from 
(iv) of Remark 1.2, $(5.3)$, (i) and (iii) of Lemma 5.2 that, for each unit normal vector $w$ of 
$F$ and each $\mu\in{\rm Spec}_JR(w)\setminus\{0\}$, 
${\rm Ker}(A^F_w\pm\sqrt{-\mu}I)\cap{\rm Ker}(R(w)-\mu I)=\{0\}$ holds.  
Therefore, it follows from Lemma 4.1 that $\widetilde F$ is 
a proper anti-Kaehlerian Fredholm submanifold and, for each unit normal 
vector $w$ of $F$, we have 
${\rm Tr}_JA^{\widetilde F}_{w^L}={\rm Tr}_JA^F_w$.  
It is clear that 
$\widetilde{\psi}_{\vert r_0\vert}(\frac{r_0}{\vert r_0\vert}v^L_{u_0})$ is 
the horizontal lift of $\psi_{\vert r_0\vert}
(\frac{r_0}{\vert r_0\vert}v_{x_0})$ to $\widetilde f_{r_0}(u_0)$.  
Hence we have 
$${\rm Tr}_JA^F_{\psi_{\vert r_0\vert}(\frac{r_0}{\vert r_0\vert}v_{x_0})}=
{\rm Tr}_JA^{\widetilde F}_{\widetilde{\psi}_{\vert r_0\vert}(\frac{r_0}
{\vert r_0\vert}v^L_{u_0}),}\leqno{(5.5)}$$, 
From $(5.4)$ and $(5.5)$, we have 
${\rm Tr}_JA^F_{\psi_{\vert r_0\vert}(\frac{r_0}{\vert r_0\vert}v_{x_0})}=0$.  
This completes the proof.  \hspace{0.7truecm}q.e.d.

\vspace{0.5truecm}

Now we prepare the following lemma to prove Theorem C.  

\vspace{0.5truecm}

\noindent
{\bf Lemma 5.3.} {\sl Let $M$ be a curvature-adapted isoparametric 
$C^{\omega}$-hypersurface in a symmetric space $N:=G/K$ of non-compact type.  
Assume that $M$ has no focal point of non-Euclidean type on $N(\infty)$.  
Then, for any complex focal radius $r$ of $M$, we have 
$${\rm Spec}\left(A_x\vert_{{\rm Ker}\,R(v_x)}\right)\subset 
\left\{\frac{1}{{\rm Re}\,r},\,\,0\right\}$$
and 
$${\rm Spec}\left(A_x\vert_{{\rm Ker}(R(v_x)-\mu I)}\right)\subset 
\left\{\frac{\sqrt{-\mu}}{\tanh(\sqrt{-\mu}{\rm Re}\,r)},\,\,
\sqrt{-\mu}\tanh(\sqrt{-\mu}{\rm Re}\,r)\right\}$$ for 
$\mu\in{\rm Spec}R(v_x)\setminus\{0\}$, 
where $x$ is an arbitrary point of $M$.}

\vspace{0.5truecm}

\noindent
{\it Proof.} 
For simplicity, we set $D_{\mu}:={\rm Ker}(R(v_x)-\mu\,{\rm id})$ for 
each $\mu\in{\rm Spec}\,R(v_x)$.  
Let $r_0$ be the complex focal radius of $M$ with 
${\rm Re}r_0=\displaystyle{\mathop{\max}_r{\rm Re}r}$, where $r$ 
runs over the set of all complex focal radii of $M$.  Let $(\lambda,\mu)
\in S_{r_0}^x\setminus\{(0,0)\}$ and $r$ a complex focal radius including 
${\rm Ker}(A_v-\lambda I)\cap D_{\mu}$ as the focal space, 
that is, $\lambda=\hat{\tau}_r(\mu)$ (see (ii) of Remark 1.2).  
Set $c_{\lambda,\mu}:=
-\frac{\mu+\lambda\hat{\tau}_{r_0}(\mu)}
{\lambda-\hat{\tau}_{r_0}(\mu)}$.  
We shall show ${\rm Re}\,c_{\lambda,\mu}\leq0$.  The argument divides into 
the following three cases:

\vspace{0.2truecm}

\centerline{${\rm(i)}\,\, \mu=0\quad\,\,
{\rm(ii)}\,\,0<\sqrt{-\mu}<\vert\lambda\vert
\quad\,\,{\rm(iii)}\,\,\vert\lambda\vert<\sqrt{-\mu}.$}

\vspace{0.2truecm}

\noindent
First we consider the case (i).  
Then we have 
$c_{\lambda,\mu}=\frac{\lambda}{1-\lambda r_0}$.  
Also, we can show $\lambda=\frac{1}{r}$.  Hence we have 
$$c_{\lambda,\mu}=\frac{1}{r-r_0}.\leqno{(5.6)}$$
Furthermore, we have ${\rm Re}\,c_{\lambda,\mu}\leq0$ from the choice of 
$r_0$.  Next we consider the case (ii).  Since $\lambda=\hat{\tau}_r(\mu)$ and 
$\lambda$ is a real number with $\vert\lambda\vert>\sqrt{-\mu}$, we can show 
$\lambda=\hat{\tau}_{{\rm Re}\,r}(\mu)
(=\frac{\sqrt{-\mu}}{\tanh(\sqrt{-\mu}{\rm Re}\,r)})$ and 
$r\equiv{\rm Re}\,r$ (${\rm mod}\,\frac{\pi{\bf i}}{\sqrt{-\mu}}$).  
Hence we have $c_{\lambda,\mu}=\hat{\tau}_{(r_0-{\rm Re}\,r)}(\mu)$, 
where we note that ${\rm Re}r\not\equiv r_0$ 
(${\rm mod}\,\frac{\pi{\bf i}}{\sqrt{-\mu}}$) because 
$(\lambda,\mu)\in S_{r_0}^x$.  Therefore, we obtain 
$${\rm Re}\,c_{\lambda,\mu}=\frac{\sqrt{-\mu}
\left(1+\tan^2(\sqrt{-\mu}{\rm Im}r_0)\right)\tanh(\sqrt{-\mu}
({\rm Re}r-{\rm Re}r_0))}{\tanh^2(\sqrt{-\mu}
({\rm Re}r-{\rm Re}r_0))+\tan^2(\sqrt{-\mu}{\rm Im}r_0)}\leq0
\leqno{(5.7)}$$
because ${\rm Re}r\leq{\rm Re}r_0$.  Next we consider the case (iii).  
Since $\lambda=\hat{\tau}_r(\mu)$ and 
$\lambda$ is a real number with $\vert\lambda\vert<\sqrt{-\mu}$, we can show 
$\lambda=\hat{\tau}_{({\rm Re}\,r+\frac{\pi{\bf i}}{2\sqrt{-\mu}})}(\mu)
(=\sqrt{-\mu}\tanh(\sqrt{-\mu}{\rm Re}\,r))$ and 
$r\equiv{\rm Re}\,r+\frac{\pi{\bf i}}{2\sqrt{-\mu}}$ 
(${\rm mod}\,\frac{\pi{\bf i}}{\sqrt{-\mu}}$).  
Hence we have $c_{\lambda,\mu}
=\hat{\tau}_{(r_0-{\rm Re}r+\frac{\pi{\bf i}}{2\sqrt{-\mu}})}(\mu)$.  
Therefore, we obtain 
$${\rm Re}c_{\lambda,\mu}=\frac{\sqrt{-\mu}\left(1+\tan^2
(\sqrt{-\mu}{\rm Im}r_0)\right)\tanh(\sqrt{-\mu}
({\rm Re}r-{\rm Re}r_0))}{1+\tanh^2(\sqrt{-\mu}({\rm Re}r-{\rm Re}r_0))
\tan^2(\sqrt{-\mu}{\rm Im}r_0)}\leq0. \leqno{(5.8)}$$
Thus ${\rm Re}c_{\lambda,\mu}\leq0$ is shown in general.  Hence, 
from the identity in Theorem B, ${\rm Re}c_{\lambda,\mu}=0$ 
($(\lambda,\mu)\in S_{r_0}^x$) follows, where we note that $c_{0,0}=0$.  
In case of (i), it follows from $(5.6)$ that 
$\displaystyle{{\rm Re}\left(\frac{1}{r-r_0}\right)=0}$.  
Hence we have ${\rm Re}\,r={\rm Re}\,r_0(<\infty)$ or $r=\infty$.  
If ${\rm Re}\,r={\rm Re}\,r_0(<\infty)$, then we have $\lambda=\frac{1}{r}
=\frac{1}{{\rm Re}\,r_0}=\hat{\tau}_{{\rm Re}\,r_0}(0)$ (which does not happen 
if $r_0$ is real because $(\lambda,0)\in S_{r_0}^x$).  
Also, if $r=\infty$, then we have $\lambda=0$.  
Thus we have 
$${\rm Spec}(A_x\vert_{D_0})\subset\left\{\frac{1}{{\rm Re}\,r_0},\,0\right\}.
\leqno{(5.9)}$$
In case of (ii), it follows from (5.7) that ${\rm Re}r={\rm Re}r_0$.  
Hence we have $\lambda=\hat{\tau}_{{\rm Re}\,r_0}(\mu)$ (which does not happen 
if $r_0\equiv{\rm Re}\,r_0$ (${\rm mod}\,\frac{\pi{\bf i}}{\sqrt{-\mu}}$) 
because $(\lambda,\mu)\in S_{r_0}^x$).  
In case of (iii), it follows from (5.8) that ${\rm Re}r={\rm Re}r_0$.  
Hence we have 
$\lambda=\hat{\tau}_{({\rm Re}\,r_0+\frac{\pi{\bf i}}{2\sqrt{-\mu}})}(\mu)$ 
(which does not happen if $r_0\equiv{\rm Re}\,r_0+
\frac{\pi{\bf i}}{2\sqrt{-\mu}}$ (${\rm mod}\,\frac{\pi{\bf i}}{\sqrt{-\mu}}$)
because $(\lambda,\mu)\in S_{r_0}^x$).  
Hence we have 
$${\rm Spec}(A_x\vert_{D_{\mu}})
\subset\left\{\frac{\sqrt{-\mu}}{\tanh(\sqrt{-\mu}{\rm Re}r_0)},\,
\sqrt{-\mu}\tanh(\sqrt{-\mu}{\rm Re}r_0)\right\}.\leqno{(5.10)}$$
This complets the proof.  
\hspace{9.75truecm}q.e.d.

\vspace{0.5truecm}

Next we prove Theorem C in terms of this Lemma and its proof.  

\vspace{0.5truecm}

\noindent
{\it Proof of Theorem C.} 
According to the proof of Lemma 5.3, the real parts of complex focal radii of 
$M$ coincide with one another.  Denote by $s_0$ this real part.  
Then, according to Lemma 5.3, we have 
$${\rm Spec}(A_x\vert_{D_0})\subset\left\{\frac{1}{s_0},0\right\}$$
and 
$${\rm Spec}(A_x\vert_{D_{\mu}})\subset
\left\{\frac{\sqrt{-\mu}}{\tanh(\sqrt{-\mu}s_0)},\,\,
\sqrt{-\mu}\tanh(\sqrt{-\mu}s_0)\right\}\,\,\,\,
(\mu\in{\rm Spec}\,R(v_x)\setminus\{0\}).$$
Set $\displaystyle{D_0^V:={\rm Ker}\left(A_x\vert_{D_0}-\frac{1}{s_0}{\rm id}
\right)}$, $D_0^H:={\rm Ker}A_x\vert_{D_0}$, 
$$D_{\mu}^V:={\rm Ker}\left(A_x\vert_{D_{\beta}}-
\frac{\sqrt{-\mu}}{\tanh(\sqrt{-\mu}s_0)}\,{\rm id}\right)$$
and 
$$D_{\mu}^H:={\rm Ker}\left(A_x\vert_{D_{\beta}}-
\sqrt{-\mu}\tanh(\sqrt{-\mu}s_0)\,{\rm id}\right).$$
According to (ii) of Remark 1.2, if 
$\displaystyle{D_0^V\oplus\left(
\mathop{\oplus}_{\mu\in{\rm Spec}\,R(v_x)\setminus\{0\}}D_{\mu}^V\right)
\not=\{0\}}$, then $s_0$ is a (real) focal radius of $M$ 
whose focal space is equal to $\displaystyle{D_0^V\oplus\left(
\mathop{\oplus}_{\mu\in{\rm Spec}\,R(v_x)\setminus\{0\}}D_{\mu}^V\right)
\not=\{0\}}$.  Let $\eta_{sv}$ ($s\in{\Bbb R}$) be the end-point map for 
$sv$.  Set $M_s:=\eta_{sv}(M)$.  
Set $F:=M_{s_0}$.  If $s_0$ is a (real) focal radius of $M$, then $F$ is 
the only focal submanifold of $M$, and if $s_0$ is not a (real) focal radius 
of $M$, then $F$ is a parallel submanifold of $M$.  
Without loss of generality, we may assume that $eK\in F$.  
Define a unit normal vector field $v^s$ of 
$M_s$ ($0\leq s<s_0$) by $v^s_{\eta_{sv}(x)}=\gamma_{v_x}'(s)$ ($x\in M$).  
Denote by $A^s$ ($0\leq s<s_0$) the shape operator of $M_s$ (for $v^s$) and 
$A^F$ the shape tensor of $F$.  
Set $(D_0^V)^s:=(\eta_{sv})_{\ast}(D_0^V)$ ($0\leq s<s_0$) and 
$(D_{\mu}^V)^s:=(\eta_{sv})_{\ast}(D_{\mu}^V)$ ($0\leq s<s_0,\,
\mu\in{\rm Spec}\,R(v_x)\setminus\{0\}$).  
Also, set $(D_0^H)^s:=(\eta_{sv})_{\ast}(D_0^H)$ ($s\in{\Bbb R}$) 
and $(D_{\mu}^H)^s:=(\eta_{sv})_{\ast}(D_{\mu}^H)$ ($s\in{\Bbb R},\,
\mu\in{\rm Spec}\,R(v_x)\setminus\{0\}$).  Easily we have 
$$T_{\eta_{s_0v}(x)}F=(D_0^H)^{s_0}_{\eta_{s_0v}(x)}\oplus
\left(\mathop{\oplus}_{\mu\in{\rm Spec}\,R(v_x)\setminus\{0\}}
(D_{\mu}^H)^{s_0}_{\eta_{s_0v}(x)}\right).\leqno{(5.11)}$$
Also, we can show 
$$
A^s_{\eta_{sv}(x)}\vert_{(D_0^H)^s_{\eta{sv}(x)}}=0\,\,\,\,(0\leq s<s_0)$$
and 
$$A^s_{\eta_{sv}(x)}\vert_{(D_{\beta}^H)^s_{\eta_{sv}(x)}}
=\mu\tanh(\sqrt{-\mu}(s_0-s))\,{\rm id}\,\,\,\,
(0\leq s<s_0) .$$
Hence we have 
$$A^F_{\psi_{s_0}(v_x)}\vert_{(D_0^H)^{s_0}_{\eta_{s_0v}(x)}}=0$$
and 
$$A^F_{\psi_{s_0}(v_x)}
\vert_{(D_{\beta}^H)^{s_0}_{\eta_{s_0v}(x)}}
=\left(\lim_{s\to s_0-0}\sqrt{-\mu}\tanh(\sqrt{-\mu}(s_0-s))\right)
{\rm id}=0,$$
where $\psi$ is the geodesic flow of $G/K$.  
From these relations and $(5.11)$, we obtain 
$A^F_{\psi_{s_0}(v_x)}=0$.  
Since this relation holds for any $x\in M$, $F$ is totally geodesic.  
Denote by $\exp^{\perp}$ the normal exponential map for $F$.  
Since the real parts of complex focal radii of $M$ coincide with one another, 
the normal umbrella $\exp^{\perp}(T^{\perp}_xF)$'s ($x\in F$) do not intersect with 
one another.  
From this fact, an involutive diffeomorphism $\tau:G/K\to G/K$ having $F$ 
as the fixed point set is well-defined by 
$\tau(\exp^{\perp}(w)):=\exp^{\perp}(-w)\,\,\,\,(w\in T^{\perp}F)$.  
For each $s\in{\Bbb R}\setminus\{s_0\}$, the restriction $\tau\vert_{M_s}$ of 
$\tau$ to $M_s$ coincides with the end-point map $\eta_{2(s_0-s)v^s}$ for 
$2(s_0-s)v^s$.  
Since $F$ is totally geodesic, we see that $\eta_{2(s_0-s)v^s}$ (hence $\tau\vert_{M_s}$) 
is an isometry of $M_s$.  From this fact, it follows that $\tau$ is an isometry of $G/K$.  
Hence $F$ is reflective.  
Furthermore, by imitating the proof of Proposition 1.12 of [KiT], we can show that 
$F$ is an orbit of a Hermann action on $G/K$ as follows.  
Take ${\rm Exp}\,Z_0\in F$, where ${\rm Exp}$ is the exponential map of $G/K$ at $o$.  
Set $\mathfrak m:={\rm Ad}(\exp(-Z_0))((\exp\,Z_0)_{\ast}^{-1}(T_{{\rm Exp}\,Z_0}F))$, 
where ${\rm Ad}$ is the adjoint operator of $G$.  
Define a subalgebra $\mathfrak k'$ of $\mathfrak g$ by 
$\mathfrak k':=\{X\in\mathfrak k\,\vert\,{\rm ad}(X)\mathfrak m=\mathfrak m\}$ and 
set $\mathfrak h:=\mathfrak k'+\mathfrak m$, which is a subalgebra of $\mathfrak g$.  
Set $H:=I(\exp\,Z_0)(\exp(\mathfrak h))$, where $I(\exp\,Z_0)$ is the inner automorphism of 
$G$ by $\exp\,Z_0$.  Easily we can show that 
$T_{{\rm Exp}\,Z_0}(H{\rm Exp}\,Z_0)=T_{{\rm Exp}\,Z_0}F$ and hence $H{\rm Exp}\,Z_0=F$.  
Define an involution $\hat{\tau}$ of $G$ by $\hat{\tau}(g):=\tau\circ g\circ\tau^{-1}$ 
($g\in G$).  It is easy to show that $({\rm Fix}\,\hat{\tau})_0\subset H\subset{\rm Fix}\,
\hat{\tau}$.  Thus $H\curvearrowright G/K$ is a Hermann action.  
Let $H^{\bf c}$ be the complexification of $H$ and $M^{\bf c}(\subset G^{\bf c}/K^{\bf c})$ 
be the complete complexification of $M$.  See [Koi6] about the definition of the complete 
complexification of $M$.  
Since both $H^{\bf c}\cdot o$ and $M^{\bf c}$ are anti-Kaehler equifocal submanifolds 
having $F^{\bf c}$ as a focal submanifold, they are equal to one of the partial tubes over 
$F^{\bf c}$ stated in Section 5 in [Koi6].  Thus they coincides with each other.  
Furthermore, from this fact, we can derive $H\cdot o=M$.  
This completes the proof.  \hspace{10.3truecm}q.e.d.


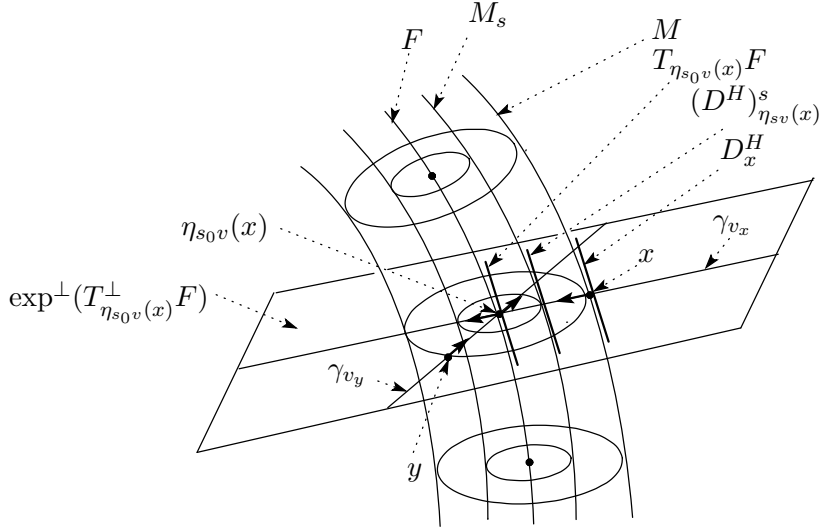
\vspace{1truecm}

\centerline{
\unitlength 0.1in
\begin{picture}( 70.3100, 34.7500)(-20.1400,-40.5000)
%
\special{pn 8}%
\special{ar 3036 1354 488 204  0.1722914 0.1999156}%
%
\special{pn 8}%
\special{ar 1940 3586 1132 2252  5.0713373 6.1832095}%
%
\special{pn 8}%
\special{pa 2196 2138}%
\special{pa 1794 2974}%
\special{fp}%
%
\special{pn 8}%
\special{pa 1794 2974}%
\special{pa 4640 2322}%
\special{fp}%
%
\special{pn 8}%
\special{pa 4640 2322}%
\special{pa 5018 1538}%
\special{fp}%
\special{pa 5018 1538}%
\special{pa 5018 1538}%
\special{fp}%
%
\special{pn 8}%
\special{pa 5018 1548}%
\special{pa 3752 1814}%
\special{fp}%
%
\special{pn 8}%
\special{pa 2208 2138}%
\special{pa 2720 2036}%
\special{fp}%
%
\special{pn 8}%
\special{pa 2792 2016}%
\special{pa 3680 1834}%
\special{fp}%
%
\special{pn 8}%
\special{pa 2792 2740}%
\special{pa 3936 1772}%
\special{fp}%
%
\special{pn 20}%
\special{sh 1}%
\special{ar 3108 2476 10 10 0  6.28318530717959E+0000}%
\special{sh 1}%
\special{ar 3108 2476 10 10 0  6.28318530717959E+0000}%
%
\special{pn 20}%
\special{sh 1}%
\special{ar 3850 2150 10 10 0  6.28318530717959E+0000}%
\special{sh 1}%
\special{ar 3850 2150 10 10 0  6.28318530717959E+0000}%
%
\special{pn 8}%
\special{pa 4130 814}%
\special{pa 3368 1204}%
\special{dt 0.045}%
\special{sh 1}%
\special{pa 3368 1204}%
\special{pa 3436 1190}%
\special{pa 3416 1180}%
\special{pa 3418 1156}%
\special{pa 3368 1204}%
\special{fp}%
%
\special{pn 8}%
\special{pa 4020 3004}%
\special{pa 4018 3036}%
\special{pa 4006 3066}%
\special{pa 3986 3092}%
\special{pa 3964 3114}%
\special{pa 3938 3134}%
\special{pa 3910 3148}%
\special{pa 3882 3164}%
\special{pa 3852 3176}%
\special{pa 3822 3188}%
\special{pa 3792 3198}%
\special{pa 3762 3206}%
\special{pa 3730 3214}%
\special{pa 3700 3222}%
\special{pa 3668 3228}%
\special{pa 3636 3232}%
\special{pa 3604 3236}%
\special{pa 3572 3240}%
\special{pa 3540 3242}%
\special{pa 3508 3244}%
\special{pa 3476 3244}%
\special{pa 3444 3244}%
\special{pa 3412 3242}%
\special{pa 3380 3238}%
\special{pa 3348 3236}%
\special{pa 3318 3232}%
\special{pa 3286 3226}%
\special{pa 3254 3218}%
\special{pa 3224 3210}%
\special{pa 3194 3200}%
\special{pa 3164 3188}%
\special{pa 3136 3174}%
\special{pa 3108 3156}%
\special{pa 3084 3136}%
\special{pa 3064 3110}%
\special{pa 3052 3080}%
\special{pa 3052 3048}%
\special{pa 3062 3018}%
\special{pa 3080 2992}%
\special{pa 3102 2968}%
\special{pa 3126 2948}%
\special{pa 3154 2932}%
\special{pa 3182 2916}%
\special{pa 3210 2902}%
\special{pa 3240 2890}%
\special{pa 3270 2880}%
\special{pa 3302 2872}%
\special{pa 3332 2864}%
\special{pa 3364 2856}%
\special{pa 3396 2850}%
\special{pa 3426 2844}%
\special{pa 3458 2840}%
\special{pa 3490 2836}%
\special{pa 3522 2834}%
\special{pa 3554 2832}%
\special{pa 3586 2832}%
\special{pa 3618 2832}%
\special{pa 3650 2834}%
\special{pa 3682 2836}%
\special{pa 3714 2840}%
\special{pa 3746 2844}%
\special{pa 3778 2848}%
\special{pa 3808 2856}%
\special{pa 3840 2864}%
\special{pa 3870 2874}%
\special{pa 3900 2884}%
\special{pa 3928 2898}%
\special{pa 3956 2916}%
\special{pa 3982 2934}%
\special{pa 4002 2958}%
\special{pa 4016 2988}%
\special{pa 4020 3004}%
\special{sp}%
%
\special{pn 8}%
\special{pa 3826 2150}%
\special{pa 3830 2182}%
\special{pa 3824 2212}%
\special{pa 3808 2240}%
\special{pa 3790 2266}%
\special{pa 3768 2290}%
\special{pa 3742 2310}%
\special{pa 3716 2328}%
\special{pa 3688 2346}%
\special{pa 3662 2362}%
\special{pa 3632 2376}%
\special{pa 3604 2388}%
\special{pa 3574 2400}%
\special{pa 3544 2412}%
\special{pa 3514 2422}%
\special{pa 3482 2432}%
\special{pa 3452 2440}%
\special{pa 3420 2446}%
\special{pa 3390 2452}%
\special{pa 3358 2458}%
\special{pa 3326 2462}%
\special{pa 3294 2466}%
\special{pa 3262 2468}%
\special{pa 3230 2470}%
\special{pa 3198 2472}%
\special{pa 3166 2472}%
\special{pa 3134 2470}%
\special{pa 3102 2466}%
\special{pa 3070 2460}%
\special{pa 3040 2454}%
\special{pa 3008 2448}%
\special{pa 2978 2436}%
\special{pa 2950 2422}%
\special{pa 2924 2404}%
\special{pa 2900 2382}%
\special{pa 2884 2354}%
\special{pa 2878 2322}%
\special{pa 2886 2292}%
\special{pa 2900 2262}%
\special{pa 2916 2236}%
\special{pa 2940 2212}%
\special{pa 2964 2194}%
\special{pa 2990 2174}%
\special{pa 3018 2156}%
\special{pa 3044 2140}%
\special{pa 3074 2126}%
\special{pa 3102 2112}%
\special{pa 3132 2100}%
\special{pa 3162 2088}%
\special{pa 3192 2078}%
\special{pa 3224 2070}%
\special{pa 3254 2062}%
\special{pa 3284 2054}%
\special{pa 3316 2048}%
\special{pa 3348 2042}%
\special{pa 3380 2038}%
\special{pa 3412 2034}%
\special{pa 3444 2030}%
\special{pa 3476 2028}%
\special{pa 3508 2028}%
\special{pa 3540 2028}%
\special{pa 3572 2030}%
\special{pa 3602 2034}%
\special{pa 3634 2038}%
\special{pa 3666 2044}%
\special{pa 3698 2050}%
\special{pa 3728 2062}%
\special{pa 3756 2076}%
\special{pa 3784 2092}%
\special{pa 3808 2114}%
\special{pa 3822 2142}%
\special{pa 3826 2150}%
\special{sp}%
%
\special{pn 8}%
\special{pa 3454 1358}%
\special{pa 3466 1388}%
\special{pa 3468 1420}%
\special{pa 3460 1452}%
\special{pa 3446 1480}%
\special{pa 3430 1508}%
\special{pa 3408 1532}%
\special{pa 3386 1556}%
\special{pa 3364 1578}%
\special{pa 3340 1598}%
\special{pa 3314 1618}%
\special{pa 3288 1636}%
\special{pa 3260 1652}%
\special{pa 3232 1668}%
\special{pa 3204 1684}%
\special{pa 3176 1698}%
\special{pa 3146 1712}%
\special{pa 3116 1724}%
\special{pa 3088 1734}%
\special{pa 3058 1746}%
\special{pa 3026 1756}%
\special{pa 2996 1764}%
\special{pa 2964 1772}%
\special{pa 2934 1778}%
\special{pa 2902 1784}%
\special{pa 2870 1790}%
\special{pa 2838 1792}%
\special{pa 2806 1794}%
\special{pa 2774 1794}%
\special{pa 2742 1792}%
\special{pa 2710 1788}%
\special{pa 2680 1782}%
\special{pa 2648 1772}%
\special{pa 2620 1756}%
\special{pa 2596 1736}%
\special{pa 2578 1710}%
\special{pa 2568 1680}%
\special{pa 2570 1648}%
\special{pa 2580 1616}%
\special{pa 2594 1588}%
\special{pa 2612 1560}%
\special{pa 2632 1536}%
\special{pa 2654 1514}%
\special{pa 2678 1492}%
\special{pa 2704 1472}%
\special{pa 2728 1454}%
\special{pa 2756 1436}%
\special{pa 2782 1418}%
\special{pa 2810 1402}%
\special{pa 2838 1386}%
\special{pa 2866 1372}%
\special{pa 2896 1360}%
\special{pa 2926 1348}%
\special{pa 2956 1336}%
\special{pa 2986 1326}%
\special{pa 3016 1318}%
\special{pa 3048 1310}%
\special{pa 3078 1302}%
\special{pa 3110 1296}%
\special{pa 3142 1290}%
\special{pa 3174 1284}%
\special{pa 3206 1282}%
\special{pa 3238 1282}%
\special{pa 3270 1282}%
\special{pa 3302 1282}%
\special{pa 3332 1288}%
\special{pa 3364 1296}%
\special{pa 3394 1308}%
\special{pa 3422 1322}%
\special{pa 3444 1346}%
\special{pa 3454 1358}%
\special{sp}%
%
\special{pn 8}%
\special{ar 2488 3810 1606 3118  5.1603701 6.1236300}%
%
\special{pn 8}%
\special{ar 1454 3340 2324 2954  5.4241057 6.2783132}%
%
\special{pn 8}%
\special{ar 1418 3320 1898 2548  5.3585647 6.2831853}%
\special{ar 1418 3320 1898 2548  0.0000000 0.0020052}%
%
\special{pn 8}%
\special{pa 3594 2200}%
\special{pa 3592 2230}%
\special{pa 3572 2256}%
\special{pa 3548 2278}%
\special{pa 3520 2294}%
\special{pa 3492 2306}%
\special{pa 3462 2320}%
\special{pa 3432 2328}%
\special{pa 3400 2336}%
\special{pa 3368 2342}%
\special{pa 3338 2346}%
\special{pa 3306 2346}%
\special{pa 3274 2346}%
\special{pa 3242 2342}%
\special{pa 3210 2336}%
\special{pa 3182 2322}%
\special{pa 3160 2298}%
\special{pa 3158 2266}%
\special{pa 3176 2240}%
\special{pa 3200 2218}%
\special{pa 3226 2202}%
\special{pa 3254 2186}%
\special{pa 3284 2176}%
\special{pa 3316 2166}%
\special{pa 3346 2158}%
\special{pa 3378 2152}%
\special{pa 3410 2148}%
\special{pa 3442 2146}%
\special{pa 3474 2146}%
\special{pa 3506 2148}%
\special{pa 3536 2156}%
\special{pa 3566 2168}%
\special{pa 3590 2188}%
\special{pa 3594 2200}%
\special{sp}%
%
\special{pn 8}%
\special{pa 3752 3004}%
\special{pa 3746 3036}%
\special{pa 3724 3058}%
\special{pa 3698 3076}%
\special{pa 3668 3088}%
\special{pa 3638 3100}%
\special{pa 3608 3108}%
\special{pa 3576 3114}%
\special{pa 3544 3120}%
\special{pa 3512 3122}%
\special{pa 3480 3122}%
\special{pa 3448 3120}%
\special{pa 3416 3118}%
\special{pa 3384 3112}%
\special{pa 3354 3100}%
\special{pa 3328 3084}%
\special{pa 3308 3058}%
\special{pa 3310 3028}%
\special{pa 3330 3002}%
\special{pa 3356 2982}%
\special{pa 3384 2970}%
\special{pa 3414 2958}%
\special{pa 3444 2950}%
\special{pa 3476 2942}%
\special{pa 3508 2938}%
\special{pa 3540 2934}%
\special{pa 3572 2934}%
\special{pa 3604 2936}%
\special{pa 3636 2938}%
\special{pa 3668 2944}%
\special{pa 3696 2954}%
\special{pa 3726 2968}%
\special{pa 3748 2992}%
\special{pa 3752 3004}%
\special{sp}%
%
\special{pn 8}%
\special{pa 3218 1436}%
\special{pa 3222 1466}%
\special{pa 3210 1496}%
\special{pa 3190 1520}%
\special{pa 3166 1542}%
\special{pa 3140 1560}%
\special{pa 3112 1576}%
\special{pa 3084 1590}%
\special{pa 3054 1604}%
\special{pa 3024 1614}%
\special{pa 2994 1624}%
\special{pa 2962 1630}%
\special{pa 2930 1634}%
\special{pa 2898 1634}%
\special{pa 2866 1632}%
\special{pa 2836 1622}%
\special{pa 2812 1600}%
\special{pa 2808 1568}%
\special{pa 2820 1538}%
\special{pa 2842 1516}%
\special{pa 2864 1494}%
\special{pa 2890 1474}%
\special{pa 2918 1458}%
\special{pa 2946 1444}%
\special{pa 2976 1432}%
\special{pa 3006 1420}%
\special{pa 3036 1412}%
\special{pa 3068 1408}%
\special{pa 3100 1402}%
\special{pa 3132 1402}%
\special{pa 3164 1406}%
\special{pa 3194 1416}%
\special{pa 3218 1436}%
\special{sp}%
%
\special{pn 8}%
\special{ar 1320 3484 2238 3026  5.4217999 6.2384701}%
%
\special{pn 20}%
\special{sh 1}%
\special{ar 3534 3026 10 10 0  6.28318530717959E+0000}%
\special{sh 1}%
\special{ar 3534 3026 10 10 0  6.28318530717959E+0000}%
%
\special{pn 20}%
\special{sh 1}%
\special{ar 3024 1526 10 10 0  6.28318530717959E+0000}%
\special{sh 1}%
\special{ar 3024 1526 10 10 0  6.28318530717959E+0000}%
%
\special{pn 20}%
\special{sh 1}%
\special{ar 3376 2250 10 10 0  6.28318530717959E+0000}%
\special{sh 1}%
\special{ar 3376 2250 10 10 0  6.28318530717959E+0000}%
%
\special{pn 8}%
\special{pa 2926 2974}%
\special{pa 3108 2494}%
\special{dt 0.045}%
\special{sh 1}%
\special{pa 3108 2494}%
\special{pa 3066 2550}%
\special{pa 3088 2544}%
\special{pa 3102 2564}%
\special{pa 3108 2494}%
\special{fp}%
%
\special{pn 8}%
\special{pa 2014 2536}%
\special{pa 4836 1936}%
\special{fp}%
\special{pa 3728 2292}%
\special{pa 3728 2292}%
\special{fp}%
%
\special{pn 8}%
\special{pa 4580 1854}%
\special{pa 4458 2006}%
\special{dt 0.045}%
\special{sh 1}%
\special{pa 4458 2006}%
\special{pa 4516 1966}%
\special{pa 4492 1964}%
\special{pa 4484 1942}%
\special{pa 4458 2006}%
\special{fp}%
%
\special{pn 8}%
\special{pa 2740 2600}%
\special{pa 2886 2652}%
\special{dt 0.045}%
\special{sh 1}%
\special{pa 2886 2652}%
\special{pa 2830 2610}%
\special{pa 2836 2634}%
\special{pa 2816 2648}%
\special{pa 2886 2652}%
\special{fp}%
%
\special{pn 13}%
\special{pa 3304 1976}%
\special{pa 3472 2516}%
\special{fp}%
%
\special{pn 13}%
\special{pa 3522 1914}%
\special{pa 3692 2454}%
\special{fp}%
%
\special{pn 13}%
\special{pa 3778 1854}%
\special{pa 3948 2394}%
\special{fp}%
%
\special{pn 8}%
\special{pa 4520 1496}%
\special{pa 3826 1986}%
\special{dt 0.045}%
\special{sh 1}%
\special{pa 3826 1986}%
\special{pa 3892 1964}%
\special{pa 3870 1956}%
\special{pa 3868 1932}%
\special{pa 3826 1986}%
\special{fp}%
%
\special{pn 8}%
\special{pa 4094 1996}%
\special{pa 3874 2128}%
\special{dt 0.045}%
\special{sh 1}%
\special{pa 3874 2128}%
\special{pa 3942 2112}%
\special{pa 3920 2100}%
\special{pa 3922 2076}%
\special{pa 3874 2128}%
\special{fp}%
\put(28.6000,-8.6400){\makebox(0,0)[lb]{$F$}}%
\put(32.1300,-7.4500){\makebox(0,0)[lb]{$M_s$}}%
\put(41.7000,-8.0800){\makebox(0,0)[lb]{$M$}}%
\put(41.0400,-19.7500){\makebox(0,0)[lb]{$x$}}%
\put(44.9400,-18.3300){\makebox(0,0)[lb]{$\gamma_{v_x}$}}%
\put(29.7300,-30.3400){\makebox(0,0)[rt]{$y$}}%
\put(26.8300,-25.1300){\makebox(0,0)[rt]{$\gamma_{v_y}$}}%
\put(45.3100,-14.7600){\makebox(0,0)[lb]{$D^H_x$}}%
\put(43.7300,-12.6200){\makebox(0,0)[lb]{$(D^H)^s_{\eta_{sv}(x)}$}}%
\put(41.7900,-10.3800){\makebox(0,0)[lb]{$T_{\eta_{s_0v}(x)}F$}}%
%
\special{pn 8}%
\special{pa 2232 1802}%
\special{pa 3352 2240}%
\special{dt 0.045}%
\special{sh 1}%
\special{pa 3352 2240}%
\special{pa 3296 2198}%
\special{pa 3302 2222}%
\special{pa 3282 2234}%
\special{pa 3352 2240}%
\special{fp}%
\put(21.7000,-17.2100){\makebox(0,0)[rt]{$\eta_{s_0v}(x)$}}%
%
\special{pn 8}%
\special{pa 2902 886}%
\special{pa 2816 1242}%
\special{dt 0.045}%
\special{sh 1}%
\special{pa 2816 1242}%
\special{pa 2850 1182}%
\special{pa 2828 1190}%
\special{pa 2812 1174}%
\special{pa 2816 1242}%
\special{fp}%
%
\special{pn 8}%
\special{pa 3206 804}%
\special{pa 3036 1170}%
\special{dt 0.045}%
\special{sh 1}%
\special{pa 3036 1170}%
\special{pa 3082 1118}%
\special{pa 3058 1122}%
\special{pa 3046 1102}%
\special{pa 3036 1170}%
\special{fp}%
%
\special{pn 8}%
\special{pa 4544 1282}%
\special{pa 3546 1936}%
\special{dt 0.045}%
\special{sh 1}%
\special{pa 3546 1936}%
\special{pa 3614 1916}%
\special{pa 3592 1906}%
\special{pa 3592 1882}%
\special{pa 3546 1936}%
\special{fp}%
%
\special{pn 8}%
\special{pa 4264 1038}%
\special{pa 3326 1996}%
\special{dt 0.045}%
\special{sh 1}%
\special{pa 3326 1996}%
\special{pa 3388 1962}%
\special{pa 3364 1958}%
\special{pa 3358 1934}%
\special{pa 3326 1996}%
\special{fp}%
%
\special{pn 13}%
\special{pa 3836 2148}%
\special{pa 3688 2180}%
\special{fp}%
\special{sh 1}%
\special{pa 3688 2180}%
\special{pa 3756 2186}%
\special{pa 3740 2170}%
\special{pa 3748 2146}%
\special{pa 3688 2180}%
\special{fp}%
%
\special{pn 13}%
\special{pa 3112 2466}%
\special{pa 3204 2386}%
\special{fp}%
\special{sh 1}%
\special{pa 3204 2386}%
\special{pa 3140 2414}%
\special{pa 3164 2422}%
\special{pa 3166 2444}%
\special{pa 3204 2386}%
\special{fp}%
%
\special{pn 13}%
\special{pa 3362 2252}%
\special{pa 3214 2284}%
\special{fp}%
\special{sh 1}%
\special{pa 3214 2284}%
\special{pa 3282 2290}%
\special{pa 3266 2272}%
\special{pa 3274 2250}%
\special{pa 3214 2284}%
\special{fp}%
%
\special{pn 13}%
\special{pa 3390 2244}%
\special{pa 3484 2164}%
\special{fp}%
\special{sh 1}%
\special{pa 3484 2164}%
\special{pa 3420 2192}%
\special{pa 3442 2200}%
\special{pa 3446 2222}%
\special{pa 3484 2164}%
\special{fp}%
%
\special{pn 8}%
\special{pa 1894 2148}%
\special{pa 2312 2322}%
\special{dt 0.045}%
\special{sh 1}%
\special{pa 2312 2322}%
\special{pa 2258 2278}%
\special{pa 2264 2302}%
\special{pa 2244 2316}%
\special{pa 2312 2322}%
\special{fp}%
\put(18.5600,-20.7700){\makebox(0,0)[rt]{$\exp^{\perp}(T^{\perp}_{\eta_{s_0v}(x)}F)$}}%
\put(20.8000,-36.5000){\makebox(0,0)[lt]{$\displaystyle{D^H_x:=(D^H_0)_x\oplus\left(\mathop{\oplus}_{\beta\in\triangle_+\vert_{{\bf R}v_x}}(D^H_{\beta})_x\right)}$}}%
\put(20.8000,-40.5000){\makebox(0,0)[lt]{$\displaystyle{(D^H)^s_{\eta_{sv}(x)}:=(D^H_0)^s_{\eta_{sv}(x)}\oplus\left(\mathop{\oplus}_{\beta\in\triangle_+\vert_{{\bf R}v_x}}(D^H_{\beta})^s_{\eta_{sv}(x)}\right)}$}}%
\end{picture}%
\hspace{9truecm}}

\vspace{2truecm}

\centerline{{\bf Fig. 3.}}

\newpage


\centerline{{\bf References}}

\vspace{0.5truecm}

{\small

\noindent
[B1] J. Berndt, 
Real hypersurfaces with constant principal curvatures 
in complex hyperbolic 

space, J. reine angew. Math. {\bf 395} (1989) 132--141.

\noindent
[B2] J. Berndt,
Real hypersurfaces in quaterionic space forms, 
J. reine angew. Math. {\bf 419} 

(1991) 9--26.

\noindent
[BV] J. Berndt and L. Vanhecke,
Curvature adapted submanifolds, 
Nihonkai Math. J. {\bf 3} (1992) 

177--185.

\noindent
[BB] J. Berndt and M. Br$\ddot u$ck,
Cohomogeneity one actions on hyperbolic spaces, 
J. reine Angew. 

Math. {\bf 541} (2001) 209--235.

\noindent
[BCO] J. Berndt, S. Console and C. Olmos,
Submanifolds and holonomy, Research Notes 
in Mathe-

matics 434, CHAPMAN $\&$ HALL/CRC Press, Boca Raton, London, New York 
Washington, 

2003.

\noindent
[BT1] J. Berndt and H. Tamaru, Homogeneous codimension one foliations on 
noncompact sym-

metric space, J. Differential Geometry {\bf 63} (2003) 1-40.

\noindent
[BT2] J. Berndt and H. Tamaru, Cohomogeneity one actions on noncompact 
symmetric spaces 

with a totally geodesic singular orbit, 
Tohoku Math. J. {\bf 56} (2004) 163-177.


\noindent
[D] J. E. D'Atri,
Certain isoparametric families of hypersurfaces in 
symmetric spaces, J. Diffe-

rential Geometry {\bf 14} (1979) 21--40.

\noindent
[E] H. Ewert,
Equifocal submanifolds in Riemannian symmetric spaces, 
Doctoral thesis.

\noindent
[G] L. Geatti, 
Complex extensions of semisimple symmetric spaces, manuscripta math. {\bf 120} 

(2006) 1-25.

\noindent
[GT] O. Goertsches and G. Thorbergsson, On the geometry of the orbits of 
Hermann actions, 

to appear.  

\noindent
[HLO] E. Heintze, X. Liu and C. Olmos, Isoparametric submanifolds and a 
Chevalley type rest-

riction theorem, Integrable systems, geometry, and topology, 151-190, 
AMS/IP Stud. Adv. 

Math. 36, Amer. Math. Soc., Providence, RI, 2006.

\noindent
[HPTT] E. Heintze, R.S. Palais, C.L. Terng and G. Thorbergsson, 
Hyperpolar actions on symme-

tric spaces, Geometry, topology and physics for Raoul Bott (ed. S. T. Yau), 
Conf. Proc. 

Lecture Notes Geom. Topology {\bf 4}, 
Internat. Press, Cambridge, MA, 1995 pp214-245.

\noindent
[H] S. Helgason, 
Differential geometry, Lie groups and symmetric spaces, 
Academic Press, New 

York, 1978.

\noindent
[HL] W. Y. Hsiang and H. B. Lawson Jr., 
Minimal submanifolds of low cohomogeneity, 
J. Diffe-

rential Geometry {\bf 5} (1971) 1--38.

\noindent
[KiT] T. Kimura and M. S. Tanaka,
Stability of certain minimal submanifolds in compact sym-

metric spaces of rank two, Differential Geom. Appl. {\bf 27} (2009) 23--33.

\noindent
[Koi1] N. Koike,
On proper Fredholm submanifolds in a Hilbert space 
arising from submanifolds 

in a symmetric space, 
Japan. J. Math. {\bf 28} (2002) 61--80.

\noindent
[Koi2] N. Koike, 
Submanifold geometries in a symmetric space of non-compact 
type and a pseudo-

Hilbert space, Kyushu J. Math. {\bf 58} (2004) 167-202.

\noindent
[Koi3] N. Koike, 
Complex equifocal submanifolds and infinite dimensional anti-
Kaehlerian isopara-

metric submanifolds, Tokyo J. Math. {\bf 28} (2005) 201-247.

\noindent
[Koi4] N. Koike, 
Actions of Hermann type and proper complex equifocal submanifolds, 
Osaka J. 

Math. {\bf 42} (2005) 599-611.

\noindent
[Koi5] N. Koike, 
A splitting theorem for proper complex equifocal submanifolds, Tohoku Math. J. 

{\bf 58} (2006) 393-417.

\noindent
[Koi6] N. Koike, The homogeneous slice theorem for the complete 
complexification of a proper 

complex equifocal submanifold, Tokyo J. Math. {\bf 33} (2010) 1-30.

\noindent
[Koi7] N. Koike, On curvature-adapted and proper complex equifocal 
submanifolds, Kyungpook 

Math. J. {\bf 50} (2010) 509-536.

\noindent
[Koi8] N. Koike, Hermann type actions on a pseudo-Riemannian symmetric 
space, Tsukuba J. 

Math. {\bf 34} (2010) 137-172.

\noindent
[Koi9] N. Koike, Examples of a complex hyperpolar action without 
singular orbit, Cubo A Math. 

{\bf 12} (2010), 127-143.

\noindent
[Koi10] N. Koike, The complexifications of pseudo-Riemannian manifolds and 
anti-Kaehler geome-

try, arXiv:math.DG/0807.1601v3.



\noindent
[Kol] A. Kollross, A Classification of hyperpolar and cohomogeneity one 
actions, Trans. Amer. 

Math. Soc. {\bf 354} (2001) 571-612.

\noindent
[Mi] R. Miyaoka, Transnormal functions on a Riemannian manifold, 
Differential Geom. Appl. 

{\bf 31} (2013) 130-139.

\noindent
[Mu] T. Murphy, Curvature-adapted submanifolds of symmetric spaces, 
Indiana Univ. Math. J. 

(to appear) (arXiv:math.DG/1102.4756v1).

\noindent
[O] O'Neill, B.
Semi-Riemannian Geometry, with Applications to 
Relativity, Academic Press, 

New York, 1983.

\noindent
[P] R. S. Palais,
Morse theory on Hilbert manifolds, 
Topology {\bf 2} (1963) 299--340.

\noindent
[PT] R. S. Palais and C. L. Terng, Critical point theory and submanifold 
geometry, Lecture Notes 

in Math. {\bf 1353}, Springer, Berlin, 1988.

\noindent
[S1] R. Sz$\ddot{{{\rm o}}}$ke,
Adapted complex structures and 
geometric quantization, 
Nagoya Math. J. {\bf 154} 

(1999) 171--183.

\noindent
[S2] R. Sz$\ddot{{{\rm o}}}$ke,
Involutive structures on the 
tangent bundle of symmetric spaces, 
Math. Ann. {\bf 319} 

(2001) 319--348.

\noindent
[S3] R. Sz$\ddot{{{\rm o}}}$ke, 
Canonical complex structures associated to connections and 
complexifications of 

Lie groups, Math. Ann. {\bf 329} (2004), 553--591.

\noindent
[Ta] Z. Tang, 
Multiplicities of equifocal hypersurfaces in symmetric spaces, 
Asian J. Math. {\bf 2} 

(1998) 181-214.

\noindent
[TT] C. L. Terng and G. Thorbergsson, 
Submanifold geometry in symmetric spaces, 
J. Differential 

Geometry {\bf 42} (1995) 665-718.

\noindent
[Wa] Q. M. Wang, Isoparametric functions on Riemannian manifolds I, 
Math. Ann. {\bf 277} (1987) 

639-646.

\noindent
[Wu] B. Wu, Isoparametric submanifolds of hyperbolic spaces, 
Trans. Amer. Math. Soc. {\bf 331} (1992) 

609--626.
}

\vspace{1truecm}

\rightline{Department of Mathematics, Faculty of Science, }
\rightline{Tokyo University of Science}
\rightline{1-3 Kagurazaka Shinjuku-ku,}
\rightline{Tokyo 162-8601, Japan}
\rightline{(e-mail: koike@ma.kagu.tus.ac.jp)}

\end{document}